\documentclass[a4paper,11pt,twoside]{article}
\usepackage{pst-fill,pst-grad}
\usepackage{textcomp}
\usepackage[titletoc]{appendix}
\usepackage[english]{babel} 
\usepackage[utf8]{inputenc}   
\usepackage{graphicx}  
\usepackage{amsmath}
\usepackage{float} 
\usepackage{fancyhdr}
\usepackage[matrix,arrow,curve]{xy}
\usepackage{pstricks} 
\usepackage{amsmath,amsfonts,verbatim,afterpage,theorem,euscript,mathrsfs,amssymb}
\usepackage{amsfonts}
\usepackage{amssymb}
\usepackage{array}
\usepackage{dsfont}
\usepackage[colorlinks=true,linkcolor=blue,citecolor=red]{hyperref}
\usepackage{authblk} 
\usepackage{color}
\newtheorem{Definition}{Definition}[section] 
\newtheorem{Proposition}{Proposition}[section]
\newtheorem{PropositionP}{Proposition}
\newtheorem{Lemme}{Lemma}[section] 

\newtheorem{TheoremeP}{Theorem}
\newtheorem{CorollaireP}{Corollary}

\def \vu{\textbf{u}}

\def \U{\textbf{u}}

\def \V{\textbf{v}}
\def \v{\textbf{v}}
\def \Rt{\mathbb{R}^{3}}
\def \R{\mathbb{R}}

\def \finpv{\hfill $\blacksquare$  \\ \newline }
\def \pv{{\bf{Proof.}}~}

\def \ds{\displaystyle}
\def \rd{\textcolor{black}}
\def \bl{\textcolor{black}} 
\title{ \bf Liouville theorems for a stationary and   non-stationary coupled system of liquid crystal flows in local Morrey spaces} 

\author[1]{ Oscar Jarr\'in \footnote{oscar.jarrin@udla.edu.ec}} 
\affil[1]{\scriptsize Direcci\'on General de Vinculaci\'on e  Investigaci\'on  (DGVI).
	Universidad de las Américas - UDLA. Calle José Queri s/n y Av. Granados. Bloque 7, Tercer Piso, Quito, Ecuador.}  

\begin{document}

\maketitle
 
\begin{abstract}
We consider here the simplified Ericksen-Leslie system on the whole space $\Rt$. This system deals  with the incompressible Navier-Stokes equations strongly  coupled with a harmonic map flow which models the dynamical behavior for nematic liquid crystals. For  both the stationary  (time independing) case and the non-stationary (time depending) case,  using the fairly general framework of a kind of local Morrey spaces,  we  obtain some \emph{a priori} conditions on the unknowns of this coupled system to prove that they vanish identically. This results are known as  Liouville-type theorems. As a bi-product,  our theorems also improve some well-known results on  Liouville-type theorems   for the particular case of classical Navier-Stokes equations.\\[5mm]
\end{abstract}
\section{Introduction} 
In this paper, we place ourselves on the whole space $\Rt$ and   we consider a  coupled system  of the  incompressible  Navier-Stokes equations  with a harmonic map flow. This system, also known as the  \emph{simplified Ericksen-Leslie system}, was proposed by H.F. Lin in \bl{\cite{Lin}} as a simplification of the general \emph{Ericksen-Leslie system} which models the hydrodynamic flow of nematic liquid crystal material \bl{\cite{Ericksen}}, \bl{\cite{Les}}. The simplified Ericksen-Leslie system, has been successful to model various dynamical behavior for nematic liquid crystals. More precisely, it provides a well macroscopic description  of the evolution of the material under the influence of fluid velocity field, and moreover, it provides  the macroscopic description of the microscopic orientation of fluid velocity  of rod-like liquid crystals. See the book \bl{\cite{Gennes}} for more details. 

\medskip

From the mathematical point of view, the simplified Ericksen-Leslie system has recently attired a lot of interest in the research community, see, \emph{e.g.}, the articles \cite{HaLiZ,LinWang,LinLiu,LinLinWang,Wang} and the references therein, where the major  challenge is due to the \emph{strong} coupled structure of this system and the presence of a super-critical non-linear term.

\medskip

In the \emph{stationary} setting, the simplified Ericksen-Leslie system is given as follows: 
\begin{equation}\label{EickLes}
\left\{ \begin{array}{ll}\vspace{2mm} 
-\Delta \U + (\U \cdot \nabla) \, \U   + \text{div} ( \nabla \otimes  \V \odot \nabla \otimes \V) + \nabla p=0, \\ \vspace{2mm} 
-\Delta \V + (\U\cdot \nabla ) \, \V   - \vert  \nabla \otimes \V  \vert^2\, \V =0, \\
\text{div}(\U)=0. \end{array} \right.  
\end{equation}
Here,  the fluid velocity $\U:\Rt \to \Rt$, and the pressure $p:  \Rt \to \mathbb{R}$ are the classical unknowns of the fluid mechanics and this system also considers a third unknown $\V  :\Rt  \to \mathbb{S}^2$ (where $\mathbb{S}^2$ denotes the unitary sphere in $\Rt$)  which is a \emph{unit vector field}  representing the macroscopic orientation of the nematic liquid crystal molecules. For the vector field $\ds{\V=(v_i)_{1\leq i \leq 3}}$, we denote   $\ds{\nabla \otimes  \V =(\partial_i v_j)_{1\leq i,j\leq 3}}$.  In the first equation of this system, the \emph{super-critical} non-linear term: $\ds{\text{div} (\nabla \otimes \V \odot \nabla \otimes \V)}$, is given as the divergence of a symmetric tensor $\ds{ \nabla \otimes \V \odot  \nabla \otimes \V}$,  where,  for $1\leq i,j \leq 3$,  its components are defined by the expression $\ds{ (\nabla \otimes  \V \odot \nabla \otimes \V)_{i,j} = \sum_{k=1}^{3} \partial_i v_k \partial_j v_k}$, 
and then, each component of the vector field $\ds{\text{div} (\nabla \otimes \V \odot \nabla \otimes \V)}$ is explicitly   given by the following expression $\ds{\left[\text{div} ( \nabla \otimes  \V \odot  \nabla \otimes \V)\right]_i= \sum_{j=1}^{3} \sum_{k=1}^{3} \partial_{j}( \partial_i v_k \, \partial_j v_k)}$. We may observe that  due to the double derivatives in this expression,  this super-critical non-linear term is actually more \emph{delicate} to treat than the classical non-linear transport term: $\ds{(\U \cdot \nabla)\, \U }$,  and this fact makes challenging the mathematical study of  (\ref{EickLes}). See, \emph{e.g.}, the works \cite{LinWang} and \cite{LinLinWang}.

\medskip

Let us introduce the Liouville-type problem for the simplified Ericksen-Leslie system. First,  we define a weak solution of the coupled system (\ref{EickLes})  as the  triplet $(\U, p, \V)$ where:  $\U \in L^{2}_{loc}(\Rt)$, $p \in \mathcal{D}^{'}(\Rt)$ and $\V \in L^{\infty}(\Rt)$  (since by the \emph{physical model} we assume $\vert \V \vert =1$)  and  $\nabla \otimes \V \in L^{2}_{loc}(\Rt)$.  Under these hypothesis all the terms in (\ref{EickLes}) are well-defined in the distributional sense. Once we have defined the weak solutions, we may remark that the triplet $\U=0$, $p=0$ and $\nabla \otimes \V=0$ (hence $\V$ is a constant unitary vector)  is a trivial weak solution of the system (\ref{EickLes}) and it is natural to ask if this solution is the unique one (modulo constants). However, it is interesting to observe that  the answer to this question is in general  \emph{negative} and we are able to exhibit an explicit  counterexample. See Appendix \ref{AppendixA} for the computations.

\medskip

Due to the non-uniqueness of the trivial solution, we are interesting to find additional some \emph{a priori} conditions  in order to assure that the trivial solution in the unique one.  This problem is commonly known as the \emph{Liouville-type problem}. To the best of our knowledge, the first  Liouville-type  result for the coupled system (\ref{EickLes}) was  recently obtained by Y. Hao, X. Liu \& X. Zhang in \cite{HaLiZ}. In this work,  the authors consider a solution   $(\U, p, \V)$ which verifies  $\nabla \otimes  \U \in L^{2}(\Rt)$ and $\nabla \otimes \V \in L^2(\Rt)$ and moreover, under the important assumption: $\ds{\U \in L^{9/2}(\Rt)}$  and  $\ds{\nabla \otimes \V \in L^{9/2}(\Rt)}$, they obtained the identities  $\U=0$, $p=0$ and $\nabla \otimes \V=0$.   These  \emph{a priori} conditions, which actually are  \emph{decaying properties} on $\U$ and $\nabla \otimes \V$ given by the $L^{9/2}-$ norm, are interesting if we compare this result with a well-known result on the Liouville problem for the  the classical stationary and incompressible Navier-Stokes equations: 
\begin{equation}\label{NS}
 -\Delta \U + (\U \cdot \nabla) \, \U   + \nabla\,p=0, \quad \text{div}(\U)=0.
\end{equation}
For these equations, a celebrated result obtained in \bl{\cite{Galdi}} by G. Galdi shows that  if $\U \in L^{9/2}(\Rt)$ then we have $\U = 0$ and $p=0$, and then, the recent result obtained in \bl{\cite{HaLiZ}} can be regarded as a generalization of Galdi's result to the more delicate setting of the coupled system (\ref{EickLes}).

\medskip

Let us recall that the Liouville problem for the stationary Navier-Stokes equations (\ref{NS}) was extensive studied in different functional settings. Galdi's result \bl{\cite{Galdi}} was  extended to  setting of the Lorentz spaces  by H. Kozono \emph{et. al.} in \bl{\cite{KzonoAL}}. Thereafter, this work was improved to a kind of  local Lorentz-type  spaces  by G. Seregin \& W. Wang in \bl{\cite{SerWang}}. Moreover, the Liouville problem for the equations (\ref{NS}) has also largely studied in the more general setting of the Morrey spaces by D. Chamorro \emph{et. al.} in \bl{\cite{ChJaLem}} and G. Seregin in \bl{\cite{Ser2015}} and \bl{\cite{Ser2016}}. For more interesting works on the Liouville problem for the stationary Navier-Stokes equations (\ref{NS}) see the articles \bl{\cite{ChaeYoneda,ChaeWolf,ChaeWeng,KochAL}} and the references therein.

\medskip

With all this information in mind, it is quite natural to improve the Galdi's-type result for the system (\ref{EickLes}) obtained in \bl{\cite{HaLiZ}} to different functional settings. Thus, the first  aim of this paper is to study  the Liouville problem for the coupled system (\ref{EickLes}) in  a fairly general functional setting.

\medskip

A kind of \emph{local Morrey spaces} (see the expression (\ref{def-local-morrey}) below for a definition) which, roughly speaking, characterize the averaged decaying properties of functions, have recently attired the attention in the study of the existence of global in time  weak solutions  for the  classical  the Navier-Stokes \bl{\cite{BradKu,PFPG1}}, and also for the coupled system of the Magneto-hydrodynamic  equations \bl{\cite{PFOJ1,PFOJ}}.  In this paper we show that the \emph{local Morrey spaces} also  give us an interesting and general setting  to solve the Liouville problem for the coupled system (\ref{EickLes}) due to the fact   these spaces contain the classic Lebesgue spaces and the more technical Lorentz and Morrey spaces. As a bi-product, since the equations (\ref{NS}) are a particular case of the system (\ref{EickLes}) (when we set $\V$ an unitary constant vector) we also improve some well-known and recent results on the Liouville problem for (\ref{NS}).

\medskip

Our methods are essentially based in $L^p-$ local  estimates of the functions $\U$ and $ \nabla \otimes \V$, and this approach also allow us to study the Liouville problem for \emph{non-stationary} case of the coupled  system (\ref{EickLes}). Thus, in the second part of this paper, we will focus on the following Cauchy problem for the simplified Ericksen-Leslie system: 
\begin{equation}\label{EickLes-non-stat}
\left\{ \begin{array}{ll}\vspace{2mm} 
\partial_t \vu -\Delta \vu + (\vu \cdot \nabla)\, \vu  + \text{div} ( \nabla \otimes  \v \odot  \nabla  \otimes  \v) + \nabla p=0, \\ \vspace{2mm} 
\partial_t \v-\Delta \v + (\vu \cdot  \nabla)\, \v   - \vert  \nabla \otimes \v  \vert^2 \, \v =0,  \quad (\vert \v(t,x)\vert=1) \\ \vspace{2mm}
\text{div}(\vu)=0,\\ 
\vu(0,\cdot)=\vu_0, \quad \v(0, \cdot)=\v_0\, \quad \text{div}(\vu_0)=0. \end{array} \right. \\ \vspace{3mm}
\end{equation} 

For the non-stationary case the Liouville-type  problem reads as follows: if we consider the initial data $\vu_0=0$ and $\nabla \otimes \v_0=0$, \emph{i.e.}, $\v_0$ is a constant vector, then we seek  if the trivial solution  $\vu=0$ and $\nabla \otimes \v=0$ (hence is $\v$ a constant vector)  in the unique one; and thus,  we are interested in studying   some \emph{a priori} conditions on  $\vu$ and $\nabla \otimes \v$ to ensure the uniqueness of the trivial solutions arising from the data $\vu_0=0$ and $\nabla \otimes \v_0=0$.

\medskip

Our general strategy to study this problem is first to look for   some \emph{a priori} conditions on the data $\vu_0$ and $\nabla \otimes \v_0$, and moreover, some   \emph{a priori} conditions on  the solutions $\vu$ and  $\nabla \otimes \v$ in order to prove that they verify a \emph{global energy inequality} (see (\ref{global-energ}) for the details). With this global energy inequality at hand, the Liouville-type problem explained above can be easily solved when $\vu_0=0$ and $\nabla \otimes \v_0=0$. \\

Before to explain this strategy more in details, we need first to overview some previous results obtained in  the particular case  for  the Cauchy problem of the   incompressible   Navier-Stokes equations: 
\begin{equation}\label{NS-non-stat}
\left\{ \begin{array}{ll}\vspace{2mm} 
\partial_t \vu -\Delta \vu + (\vu \cdot \nabla)\, \vu   + \nabla p=0, \quad  
\text{div}(\vu)=0. \\ 
\vu(0,\cdot)=\vu_0,  \quad \text{div}(\vu_0)=0. \end{array} \right. 
\end{equation}

\medskip

For this system, J. Serrin proved in \cite{Serrin} that if $\vu_0\in L^2(\Rt)$ and if  a  Leray weak solution $\vu$ satisfies the condition  $\vu \in L^p(0,T, L
^r(\Rt))$, for  $p>2$ and $r>3$ such that $\ds{2/p + 3/r \leq 1}$, then $\vu$ verifies the \emph{global energy equality}. Thereafter, this result was generalized by H. Kozono \emph{et. al.}   in  \cite{KzonoAL}  as follows. Recall first that the notion of  \emph{weak suitable solutions} for the equations (\ref{NS-non-stat}) were introduced   in the celebrated  Cafarelli, Konh and Niremberg theory \cite{CKN}.  Then, in \cite{KzonoAL},  it is introduced the notion of \emph{generalized weak suitable solution}  (see Definition $3.1$, page $5$ of  \cite{KzonoAL}).  This notion of generalized weak suitable solution is  a generalization of the well-known  weak suitable solutions and the main difference is that it assumes  neither finite energy: $\ds{\sup_{0\leq t \leq T}\Vert \vu(t,\cdot)\Vert^{2}_{L^2}<+\infty}$, nor finite dissipation $\ds{\int_{0}^{T}\Vert \vu(t,\cdot)\Vert^{2}_{\dot{H}^1} dt <+\infty}$. In the setting of the \emph{generalized weak suitable solution},  H. Kozono \emph{et. al.}  gave a  \emph{new a priori condition} which ensures that the well-know \emph{global  energy inequality} holds. More precisely, assuming that $\vu_0 \in L^2(\Rt)$ and moreover,  within the  general framework of the Lorentz spaces and  for the  parameters $3\leq p_1, r_1, p_2, r_2\leq +\infty$ satisfying some technical conditions related to the well-known  scaling properties of the equations (\ref{NS-non-stat}), the condition  $\vu \in L^3(0,T,L^{p_1, r_1}(\Rt))\cap L^2(0,T, L^{p_2, r_2}(\Rt))$  ensures that $\vu \in L^{\infty}_{t}L^{2}_{x}\cap L^{2}_{t}\dot{H}^{1}_{x}(]0,T[\times \Rt) $ and moreover it verifies the global energy inequality, \emph{i.e.}, $\vu$ becomes a Leray weak solution.

\medskip

Following these ideas, in  Definition \ref{def-suitable} below,  we will introduce first  a notion of \emph{generalized weak suitable solutions} $(\vu, p, \v)$  for the coupled  system (\ref{EickLes-non-stat}).  Thereafter, assuming that $\vu_0 \in L^2(\Rt)$, $\v_0 \in \dot{H}^1(\Rt)$, and moreover,  using a time-space version   of the \emph{local Morrey spaces} (see the expression (\ref{def-time-space-local-Morrey}) below for a definition), we will give some \emph{a priori} conditions on $\vu$ and $\nabla \otimes \v$ to ensure that,  for a time $0<T<+\infty$ arbitrary large,  the \emph{generalized weak suitable solutions} of (\ref{EickLes-non-stat}) verify  a \emph{global energy inequality} (\ref{global-energ}). As an interesting application,  we obtain some \emph{Liouville-type results} for the non-stationary system (\ref{EickLes-non-stat}). More precisely, using the global energy inequality  we are able to prove the uniqueness of the trivial solutions $\vu=0$, $p=0$ and $ \nabla \otimes  \v=0$ for the initial data $\vu_0=0$ and $ \nabla \otimes \v_0=0$.

\medskip

This paper is organized as follows. In Section \ref{Sec:Results} below we expose all the results obtained. Then, in Section \ref{Sec:Morrey} we summarize some previous  results on the local Morrey spaces we shall use here. Section \ref{Sec:Presion} is devoted to a  characterization of the pressure term in the coupled systems (\ref{EickLes}) and (\ref{EickLes-non-stat}) which  will be useful for the next sections. Finally, in Section \ref{Sec:Stat} we study the Liouville problem for the stationary system  (\ref{EickLes}) and in Section \ref{Sec:Non-Stat} we study the Liouville problem for the non-stationary system (\ref{EickLes-non-stat}).
\section{Framework and statement of the results}\label{Sec:Results}
\subsection{The stationary case}
Recall first that in \bl{\cite{HaLiZ}}, in order to solve the Liouville problem for  (\ref{EickLes}), the authors need the additional  hypothesis on the function $\V$: $\nabla\otimes \V \in L^2(\Rt)$. In our results  we will \emph{relax} this hypothesis as follows:  for $R\geq 1$ we denote  $\mathcal{C}(R/ 2, R) = \{ x \in \Rt: R/2 < \vert x \vert < R \}$;  and from now on we will assume  
\begin{equation}\label{cond-derV}
\sup_{R\geq 1} \int_{\mathcal{C}(R/2, R)} \vert \nabla \otimes \V \vert^2 dx <+\infty. 
\end{equation}

Before to state our results, we recall the definition of the Morrey spaces and local Morrey spaces. For more references about this spaces see, \emph{e.g.}, the Chapter $8$ of the book \bl{\cite{PGLR1}} and the Section $7$ of the paper \cite{PFPG} respectively. Let  $1<p<r<+\infty$, the homogeneous Morrey space $\dot{M}^{p,r}(\Rt)$ is the set of functions $f \in L^{p}_{loc}(\Rt)$ such that 
\begin{equation}\label{Def-Morrey}
\Vert f \Vert_{\dot{M}^{p,r}}= \sup_{R>0,\,\, x_0 \in \Rt}  R^{\frac{3}{r}} \left(  \frac{1}{R^{3}} \int_{B(x_0, R)} \vert f (x) \vert^p dx \right)^{\frac{1}{p}}  < +\infty,
\end{equation} where $B(x_0,R)$ denotes the ball centered at $x_0$ and with radio $R$.  This is a homogeneous space of degree $-\frac{3}{r}$ and moreover we have the following chain of continuous embedding $L^{r}(\Rt)\subset L^{r, q}(\Rt)  \subset \dot{M}^{p,r}(\Rt)$, where, for $r \leq q \leq +\infty$ the space $L^{r,q}(\Rt)$ is a Lorentz space \bl{\cite{DCh}}.

\medskip

Observe that in  expression (\ref{Def-Morrey}) we consider the average in terms of the $L^p-$ norm of the function $f$  on the ball $B(x_0,R)$; and the term $R^{\frac{3}{q}}$ describes the decaying of this averaged quantity when $R$ is large.

\medskip

The local Morrey spaces we shall consider here describes the averaged decaying of functions in a more general setting. For $\gamma\geq 0$ and $1<p<+\infty$, we define the \emph{local Morrey space} $M^{p}_{\gamma}(\Rt)$ as the Banach space of functions $f \in L^{p}_{loc}(\Rt)$ such that 
\begin{equation}\label{def-local-morrey}
\Vert f \Vert_{M^{p}_{\gamma}} = \sup_{R \geq 1} \left( \frac{1}{R^{\gamma}} \int_{B(0,R)} \vert f (x) \vert^p dx\right)^{1/p} <+\infty. 
\end{equation}
Here the parameter $\gamma \geq 0$ characterizes the behavior of the quantity $\ds{\left(  \int_{B(0,R)} \vert f (x) \vert^p dx\right)^{1/p}}$ when $R$ is large. Moreover, for $\gamma_1 \leq \gamma_2$ we have the continuous embedding  $M^{p}_{\gamma_1}(\Rt) \subset M^{p}_{\gamma_2} (\Rt)$. Remark also that for $1 < p < r <+\infty$, setting the parameter $\gamma$ such that  $ 3( 1-p/r)<\gamma$, then we have   $\ds{\dot{M}^{p,r}}(\Rt) = M^{p}_{3( 1-p/r)}(\Rt)  \subset M^{p}_{\gamma}(\Rt)$, and in this sense the local Morrey space $M^{p}_{\gamma}(\Rt)$ is  as a generalization of the homogeneous Morrey space $\dot{M}^{p,r}(\Rt)$. 

\medskip

Finally,  we define the space $M^{p}_{\gamma,0}(\Rt)$ as the set of functions $f \in M^{p}_{\gamma}(\Rt)$ such that   
 \begin{equation}\label{cond-dec1}
 \lim_{R \to +\infty}  \left( \frac{1}{R^\gamma} \int_{\mathcal{C}(R/2, R)} \vert f (x) \vert^p dx \right)^{1/p} =0.  
 \end{equation}

\medskip

In the setting of the local Morrey spaces $M^{p}_{\gamma,0}(\Rt)$ and $M^{p}_{\gamma}(\Rt)$ defined above we set, from now on, the parameters $0<\gamma<3 \leq p <+\infty$.  The condition $0 < \gamma <3$ is required to use some useful properties of the spaces $M^{p}_{\gamma}(\Rt)$ which we summarize in Section \ref{Sec:Morrey}. On the other hand, since all our results are based on a \emph{Cacciopoly} type estimate on the term $\vec{\nabla} \otimes \U$ (for more details see the Proposition \ref{Prop-Cacciopoli}) we also need the condition $3\leq p <+\infty$. 

\medskip  

\begin{tabular}{ll}
\hspace{-1cm} 	
\begin{minipage}{10cm} 
	
With the parameters $0<\gamma<3 \leq p <+\infty$ above, we introduce now the following quantity:
\begin{equation}\label{eta}
	\eta=\eta(\gamma,p)= \frac{\gamma}{p}-\frac{3}{p}+\frac{2}{3},
\end{equation} which relates the  decaying parameter $\gamma$ with  the local integrability parameter $p$ in the definition of the spaces $M^{p}_{\gamma}(\Rt)$ given in (\ref{def-local-morrey}).  Our results deeply depends on the sign of the function $\eta(\gamma,p)$. More precisely, within the \emph{rectangular region} $(\gamma,p) \in ]0,3[\times [3,+\infty[$, we will consider first  the region  of parameters $(\gamma,p)$ where $\eta(\gamma,p) \leq 0$, and  then,  the region  of parameters $(\gamma,p)$ where $\eta(\gamma,p) > 0$.  

\medskip

In this figure, we draw these  regions. In the horizontal (red) axis we have the parameter $\gamma$, while in the vertical (green) axis we have the parameter $p$. Thus, we have $\eta(\gamma,p)>0$ in the blue sky region, we have $\eta(\gamma,p)\leq 0$ in the dark gray region, while $\eta(\gamma,p)$ is not defined in the light gray region.
	\end{minipage}\hspace{1.5cm}	
\begin{minipage}{4cm}
	\begin{center}
		\includegraphics[scale=0.38]{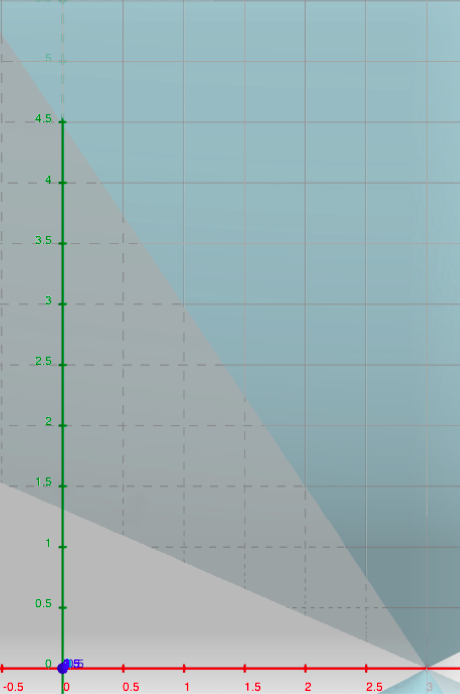} 	
	\end{center}
{\footnotesize Figure 1: Signs of  $\eta(\gamma,p)$.}
\end{minipage}	
\end{tabular}
\vspace{7mm} 

\begin{TheoremeP}\label{Th1}  Let $(\U, p, \V)$ be a smooth solution of stationary coupled system (\ref{EickLes}) such that $\V$ verifies (\ref{cond-derV}).  For  $0<\gamma<3 \leq p <+\infty$,  we assume  $\U \in M^{p}_{\gamma,0}(\Rt)$  and  $\nabla \otimes \V \in  M^{p}_{\gamma}(\Rt)$.
\begin{enumerate} 
    \item[$1)$]  If $(\gamma,p)$ are such that $\eta(\gamma,p)\leq 0$, then we have  $\U=0$,  $\nabla \otimes \V=0$ and $p$ is a constant vector.
  \item[$2)$]  If $(\gamma,p)$ are such that $\eta(\gamma,p)> 0$, and if the velocity $\U$ also verifies the following additional decaying condition:
  \begin{equation}\label{cond-dec}
  \lim_{R \to +\infty}  R^{3 \eta(\gamma,p)}\left( \frac{1}{R^{\gamma}} \int_{\mathcal{C}(R/2, R)} \vert \U(x) \vert^p dx \right)^{1/p} =0,  
  \end{equation}
  then we have  $\U=0$,  $\nabla \otimes \V=0$ and $p$ is a constant vector.
\end{enumerate} 
\end{TheoremeP}	

We observe first that among the two unknowns of the coupled system  (\ref{EickLes}), the velocity $\U$ must have a faster decaying than the derivatives of the vector field $\V$, since we have $\U \in M^{p}_{\gamma,0}(\Rt)$  and  $\nabla \otimes \V \in  M^{p}_{\gamma}(\Rt)$. 

\medskip

 In point $1)$,   for $(\gamma,p)$ such that $\eta(\gamma,p) \leq  0$, we solve the Liouville problem for the equations (\ref{EickLes}) in the  local Morrey spaces  $M^{p}_{\gamma,0}(\Rt)$ and $M^{p}_{\gamma}(\Rt)$. The main interest of this result bases on the fact that we use  a  \emph{fairly general framework} to solve this problem. Indeed, we have the following remarks.

\begin{enumerate}
\item[$\bullet$] First, if $\eta(\gamma,p)\leq 0$, then  for  $3<r<9/2$ and $3 \leq p <r$    we have the large chain of strict embedding

\[  L^r(\Rt)\subset L^{r,\infty}(\Rt)\subset \dot{M}^{p,r}(\Rt)\subset M^{p}_{\delta}(\Rt)\subset M^{p}_{\gamma,0}(\Rt),\]
involving the Lebesgue, Lorentz, Morrey and local Morrey spaces. Here, the last  embedding is due to point $1$ of Lemma \ref{Lema-Tech1} below, where we the parameter $\delta$ verifies $3(1-p/r)<\delta < \gamma$. 

\item[$\bullet$] Moreover, for the particular values $\gamma=1$ and $p=3$,  hence we have $\eta(1,3)=0$, and  for $r=9/2$ and $9/2<q<+\infty$,  we also have the embedding 
\begin{equation}\label{Emb}
L^{9/2}(\Rt)\subset L^{9/2, q}(\Rt)\subset M^{3}_{1,0}(\Rt).
\end{equation}
Indeed, if $f \in L^{9/2,q}(\Rt)$ then we have $f\in \dot{M}^{3,9/2}(\Rt)$, but due to  the identity $\dot{M}^{3,9/2}(\Rt)= M^{3}_{1}(\Rt)$,  we get $f \in M^{3}_{1}(\Rt)$. Moreover, by the following estimate: 
\[  \int_{\mathcal{C}(R/2, R)} \vert f \vert^3 dx= \int_{B(0,R)} \left\vert \mathds{1}_{\mathcal{C}(R/2,R)} f \right\vert^3 dx \leq c \, R \left\Vert  \mathds{1}_{\mathbb{C}(R/2,R)} f  \right\Vert^{3}_{L^{9/2,\infty}} \leq c \, R \left\Vert  \mathds{1}_{\mathbb{C}(R/2,R)} f  \right\Vert^{3}_{L^{9/2,q}},\]
and using the dominated convergence theorem (which is valid in the space $L^{9/2,q}(\Rt)$ for the values $9/2\leq q <+\infty$, see \cite{DCh}) we obtain: $\ds{\lim_{R \to +\infty} \frac{1}{R} \int_{\mathbb{C}(R/2,R)} \vert f \vert^3\,dx=0}$, hence we have $f \in M^{3}_{1,0}(\Rt)$. 
Due to the  embedding given in (\ref{Emb}), we may see that the recent result obtained  in \bl{\cite{HaLiZ}} for the coupled system (\ref{EickLes}) follows from point $1)$ in Theorem \ref{Th1}.

\item[$\bullet$]   Finally,   for the values in the threshold $\eta(\gamma,p)=0$,   by the expression (\ref{eta}) we have the identity  $\gamma=3- 2p/3$, and then, for $0 < \gamma < 3$ we get $3 \leq p < 9/2$.  Thus, for these values of the parameter $p$,   by the first point in Lemma \ref{Lema-Tech1} below  we  have the embedding 

\begin{equation}\label{Emb-weight}
L^{p}_{w_\gamma}(\Rt) \subset M^{p}_{\gamma,0}(\Rt),
\end{equation}
 where, for  $\ds{w_\gamma(x)= \frac{1}{(1+\vert x \vert)^\gamma}}$,  the weighted Lebesgue space $ L^{p}_{w_\gamma}(\Rt) $ is defined as  $L^{p}_{w_\delta}(\Rt)=L^{p}(w_\delta\, dx)$. 
\end{enumerate}	

In point $2)$ above,  we may observe  now that for the values $(\gamma,p)$ where $\eta(\gamma,p)>0$, the decaying properties given by  the space $M^{p}_{\gamma,0}(\Rt)$ seems not to be sufficient to solve the Liouville problem.  More precisely,  if we only consider the information $\U \in M^{p}_{\gamma,0}(\Rt)$, then the velocity $\U$ seems not to  decay at infinity fast enough and we need  to improve its decaying properties with the expression   $R^{3\, \eta(\gamma,p)}$ in (\ref{cond-dec}).  Moreover, we remark that this improvement on the decay properties (when $\eta(\gamma,p)>0$) are only needed for the velocity $\U$ and not for the function $\nabla \otimes \V$. 

\subsubsection{Some new results for the Stationary Navier-Stokes equations}

We observe that the stationary  coupled Ericksen-Leslies  system (\ref{EickLes}) contains as a particular case the stationary Navier-Stokes  equations (\ref{NS}) when setting the unitary  vector field $\V$ as a constant vector. Thus, a direct consequence of Theorem \ref{Th1} is the following new result for the equations (\ref{NS}).

\begin{CorollaireP}\label{Corollary-NS}  Let $(\U, p)$ be a smooth solution of  the stationary Navier-Stoke equations (\ref{NS}).  For  $0<\gamma<3 \leq p <+\infty$,  we assume  $\U \in M^{p}_{\gamma,0}(\Rt)$. If  $(\gamma,p)$ are such that $\eta(\gamma,p)\leq 0$, then we have  $\U=0$ and $p$ is a constant vector.
\end{CorollaireP}

It is worth mention  how this corollary improves some previous results  obtained for the Liouville problem for the stationary Navier-Stoke equations.  We observe first that  by the  embedding (\ref{Emb}) our result  improves the classical  Galdi's result  \bl{\cite{Galdi}} given in the framework of the Lebesgue space $L^{9/2}(\Rt)$  and some recent results  \bl{\cite{Jarrin}} obtained  in the framework of the Lorentz spaces $L^{9/2,q}(\Rt)$, with $9/2<q<+\infty$.  Moreover, due to the embedding $M^{p,r}(\Rt)\subset M^{p}_{\gamma,0}(\Rt)$,  with $3<r<9/2$, our result   improves some previous results obtained in the setting of the Morrey spaces in \cite{ChJaLem} and \cite{Jarrin}. Finally, due  to the embedding (\ref{Emb-weight}), our  result also improves some results proven in \cite{Phan} (see Remark $4.9$, page $10$) in the setting of weighted spaces. 
 
 \medskip 

On the other hand, we are also interested in studying the effects of removing the condition (\ref{cond-dec1}) on the velocity $\vu \in M^{p}_{\gamma,0}(\Rt)$; and we consider now  $\vu \in M^{p}_{\gamma}(\Rt)$. Within the framework of the larger space $M^{p}_{\gamma}(\Rt)$, we have the following result. 

\begin{PropositionP}\label{Prop:NS}  Let $(\U, p)$ be a smooth solution of  the stationary Navier-Stoke equations (\ref{NS}).  For  $0<\gamma<3 \leq p <+\infty$,  we assume  $\U \in M^{p}_{\gamma}(\Rt)$, and moreover, we assume that  $(\gamma,p)$ are such that $\eta(\gamma,p)\leq 0$. 
\begin{enumerate}
\item[1)] In the case when $\eta(\gamma,p)<0$, we have 	$\U=0$ and $p$ is a constant vector.
\item[2)] In the case when $\eta(\gamma,p)=0$,   if there holds  $\vu \in \dot{B}^{-1}_{\infty,\infty}(\Rt)$,  then we have  $\U=0$ and $p$ is a constant vector. 
\end{enumerate}		 
\end{PropositionP}	

We observe here that, on the one hand,  when $\eta(\gamma,p)<0$ the condition   (\ref{cond-dec1}) actually  is not required to solve the Liouville problem. On the other hand,  when $\eta(\gamma,p)=0$  we need a supplementary hypothesis on this vector field to ensure its vanishing; and this fact suggests the acuteness  of (\ref{cond-dec1}). 

 \medskip 

The supplementary condition is given in the framework of the homogeneous Besov space $\dot{B}^{-1}_{\infty, \infty}(\Rt)$, defined as the set of distributions $f \in \mathcal{S}^{'}(\Rt)$ such that $\Vert f \Vert_{\dot{B}^{-1}_{\infty, \infty}}=\sup_{t>0} t^{1/2} \Vert h_t \ast f \Vert_{L^{\infty}}<+\infty$, where $h_t$ denotes the heat kernel. This space plays a very important role in the analysis on the Navier-Stokes equations (stationary and non stationary) since this is the largest space which is invariant under scaling properties of these equations. See, for instance, the article \cite{Bourgain} and the books \cite{PGLR0} and \cite{PGLR1}.

\subsection{The non-stationary case}
From now on, let us fix a time $0<T<+\infty$.  We start by introducing the notion of \emph{generalized weak suitable solution} for the non-stationary Ericksen-Leslie system (\ref{EickLes-non-stat}). 

\begin{Definition}\label{def-suitable} Let $\vu_0 \in L^2(\Rt)$ such that $div(\vu_0)=0$ and let $\v_0 \in \dot{H}^{1}(\Rt)$.  We say that the triplet $(\vu, p, \v)$ is a generalized weak suitable solution of the coupled system  (\ref{EickLes-non-stat}) if: 
\begin{enumerate}
\item[$1)$] $\ds{\vu \in L^{3}_{loc}([0,T[\times \Rt)}$, $\ds{\nabla \otimes \vu \in L^{3}_{loc}([0,T[\times \Rt)}$ and $p \in L^{3/2}_{loc}([0,T[\times \Rt)$.
\item[$2)$] $\v \in L^{\infty}_{loc}([0,T[,L^{\infty}(\Rt))$, $\ds{\nabla \otimes \v \in L^{3}_{loc}([0,T[\times \Rt)}$ and $\Delta \v \in L^{3}_{loc}([0,T[\times \Rt)$.
\item[$3)$] The triplet $(\vu, p, \v)$ verifies the first three equations of  (\ref{EickLes-non-stat}) in $\mathcal{D}^{'}(]0,T[\times \Rt)$.
\item[$4)$] For every compact set $K \subset \Rt$, the function $\vu(t,\cdot)$ is continuous for $t \in ]0,T[$ in the weak topology of $L^2(K)$ and strongly continuous at $t=0$. Moreover, the function $\v(t,\cdot)$ is continuous for $t \in ]0,T[$ in the weak topology of $\dot{H}^{1}(K)$ and strongly continuous at $t=0$.
\item[$5)$] The triplet $(\vu, p, \v)$ verifies the following \emph{local energy inequality}:  there exist a \emph{non-negative}, locally finite measure $\mu$ on $]0,T[\times \Rt$ such that: 
\begin{equation}\label{energ-loc}
\begin{split}
\partial_t \left( \frac{\vert \vu \vert^2 +\vert \nabla \otimes \v \vert^2}{2}\right) + \vert \nabla \otimes \vu \vert^2= - \vert \Delta \v \vert^2  + \Delta \left(  \frac{\vert \vu \vert^2 +\vert \nabla \otimes \v \vert^2}{2} \right)\\
-\text{div} \left( \left[ \frac{\vert \vu \vert^2 +\vert \nabla \otimes \v \vert^2}{2}+p \right] \vu \right) - \sum_{k=1}^{3} \partial_k ([\vu \cdot \nabla)\, \v ] \cdot \partial_k \v ) - \vert \nabla \otimes \v \vert^2 \v \cdot \Delta \v - \mu.
\end{split}
\end{equation}
\end{enumerate}
\end{Definition}
Observe that in point $2)$ we assume $\v \in L^{\infty}_{loc}([0,T[,L^{\infty}(\Rt))$ due to the fact that by the \emph{physical model} we have $\vert \v(t,x)\vert=1$. Observe moreover that by the hypothesis given in points $1)$ and $2)$ we have that  $\mu$ is well-defined in the distributional sense. However, the most important fact in this definition is the  positivity assumed on $\mu$ which is the whole point in the notion of \emph{suitable} solutions. 

 \medskip 

This notion of generalized weak suitable solution is close to the definition of a weak suitable solution for the coupled system (\ref{EickLes-non-stat}) given in \cite{LinLiu} (for the case of a bounded and smooth domain $\Omega \subset \Rt$). In comparison with \cite{LinLiu}, it is worth to remark that here we suppose neither $\vu \in L^{\infty}_{t}L^{2}_{x}\cap L^{2}_{t}\dot{H}^{1}_{x}$, nor $\v \in L^{\infty}_{t}\dot{H}^{1}_{x}$, and  we consider here only  \emph{locally integrable properties}. Moreover,  we assume on the pressure term a local $L^{3/2}-$ integrability, while the authors in \cite{LinLiu} assume a $L^{5/3}-$ integrability.

 \medskip 

We introduce now a time-space version of the \emph{local Morrey spaces}, for more references on these spaces see always the Section $7$ of \cite{PFPG}. For $\gamma>0$ and  $1<p<+\infty$, we define the space $M^{p}_{\gamma}L^p(0,T)$ as the Banach space of functions $f \in L^{p}_{loc}([0,T]\times \Rt)$ such that 
\begin{equation}\label{def-time-space-local-Morrey}
\Vert f \Vert_{M^{p}_{\gamma}L^p(0,T)}= \sup_{R \geq 1} \left( \frac{1}{R^\gamma} \int_{0}^{T} \int_{B(0,R)}\vert f(t,x)\vert^p dx\, dt \right)^{1/p} <+\infty.    
\end{equation}
Moreover, we define the space $\ds{M^{p}_{\gamma,0}L^p(0,T)}$ as the set of functions $f \in M^{p}_{\gamma}L^p(0,T)$ which verifies 
\begin{equation}\label{cond-dec-2}
\lim_{R \to +\infty}  \left( \frac{1}{R^\gamma} \int_{0}^{T} \int_{\mathcal{C}(R/2,R)}\vert f(t,x)\vert^p dx\, dt \right)^{1/p}=0.   
\end{equation}

\medskip

In Definition \ref{def-suitable} we observe that we need to handle the pressure $p$ and for this, before to state our next result, it is useful to give first the following characterization of the pressure term. 

\begin{Proposition}\label{Prop-presion}  Let  $(\vu, p,  \v)$ be a solution of the coupled system (\ref{EickLes-non-stat}) such that, for $0<\gamma<3$, $2<p<+\infty$, it verifies  $\vu \in M^{p}_{\gamma}L^p(0,T)$, $p \in \mathcal{D}^{'}([0,T]\times \Rt)$ and $\nabla \otimes \v \in M^{p}_{\gamma}L^p(0,T)$. Then, the  term $\nabla p$ is necessary related to $\vu$ and $\nabla\otimes \v$ through the Riesz transforms $\mathcal{R}_i= \frac{\partial_i}{\sqrt{-\Delta}}$ by the formula 
	\begin{equation}\label{Grad-Presion}
	\nabla p = \nabla\left( \sum_{i,j=1}^{3} \mathcal{R}_i \mathcal{R}_j (u_i\, u_j)+ \sum_{i,j,k=1}^{3} \mathcal{R}_i \mathcal{R}_j\left( \partial_i v_k\, \partial_j v_k\right) \right).
	\end{equation}	
\end{Proposition}	

Here, in the general setting of the time-space local Morrey spaces, we show that  the pressure  $p$ is always related to the velocity $\vu$ and the derivatives of the vector field $\V$. This results has also an independent interest when seeking for very general frameworks in which the pressure is related to the other unknowns in the equations (\ref{EickLes-non-stat}). See, for instance, \cite{BAPF} and \cite{PFPG}, for related works in the case of the Navier-Stokes equations (\ref{NS-non-stat}). 

\medskip

As mentioned in the introduction, in our next result we give some \emph{a priori} conditions on the generalized weak suitable solutions defined above to ensure that these solutions verify a global energy inequality.  

\begin{TheoremeP}\label{Th2} Let  $\vu_0 \in L^2(\Rt)$, with $div(\vu_0)=0$, and let $\v_0 \in \dot{H}^{1}(\Rt)$ be the initial data. Let   $0<T<+\infty$, and let $(\vu,p,\v)$ be a \emph{generalized weak suitable solution} of the non-stationary coupled system (\ref{EickLes-non-stat}) given in Definition \ref{def-suitable}. 
	
\medskip	
	
For $0<\gamma<3\leq p <+\infty$,   we assume  $\vu \in M^{p}_{\gamma,0}L^p(0,T)$ and $\nabla \otimes \v \in M^{p}_{\gamma}L^p(0,T)$. If $(\gamma,p)$ are such that  the quantity $\eta(\gamma,p)$  in (\ref{eta}) verifies  $\eta(\gamma,p) \leq 0$,  then we have $\vu \in L^{\infty}_{t}L^{2}_{x}\cap L^{2}_{t}\dot{H}^{1}_{x}([0,T]\times \Rt)$, $\v \in L^{\infty}_{t}\dot{H}^{1}_{x}([0,T]\times \Rt)$, and moreover, for all $t\in [0,T]$ the global energy inequality is verified: 
\begin{equation}\label{global-energ}
\Vert \vu(t,\cdot)\Vert^{2}_{L^2}+ 2 \int_{0}^{t}\Vert \vu(s,\cdot)\Vert^{2}_{\dot{H}^{1}} ds + \Vert \v(t, \cdot)\Vert^{2}_{\dot{H}^{1}} \leq \Vert \vu_0 \Vert^{2}_{L^2}+\Vert \v_0\Vert^{2}_{\dot{H}^{1}}.       
\end{equation} 
\end{TheoremeP}

As a direct application of the global energy inequality above we have the following Liouville-type result. 

\begin{CorollaireP}\label{Col1} Within the framework of Theorem \ref{Th2},  let $(\vu,p,\v)$ be a \emph{generalized weak suitable solution} of the non-stationary coupled system (\ref{EickLes-non-stat}) given by Definition \ref{def-suitable}. Moreover, let $\vu_0$ and $\v_0$ be the initial data.  If  $\vu_0=0$ and $\v_0$ is  a constant vector field, then   we  have  $\vu=0$ and $\nabla \otimes \v =0$ on $[0,T]\times \Rt$. 
\end{CorollaireP}	

To close this section, let us make the following comments. As for the stationary case, we may observe that if we set $\v_0$ and $\v$ two constant unitary vectors, then Theorem \ref{Th2} and Corollary \ref{Col1}  hold true for the classical Navier-Stokes equations (\ref{NS-non-stat}) provided that $(\vu,p)$ is a generalized weak suitable solution in the sense of Definition \ref{def-suitable} (with $\nabla \otimes \v=0$) and  $\vu \in M^{p}_{\gamma,0}L^p(0,T)$, with $(\gamma,p)$ such that $\eta(\gamma,p)\leq 0$.  In this setting, it is interesting to observe that the space $M^{p}_{\gamma,0}L^p(0,T)$  generalizes  some 
spaces in which these kind of results have been obtained in previous works \cite{Serrin}. More precisely, for  $3 < p,r \leq 9/2$   such that $2/p+3/r \leq 1$,  we have the following chain of embedding 
\begin{equation}\label{embedding}
L^{p}\Big(0,T,L^r(\Rt)\Big) \subset L^p \Big(0,T,L^{r,q}(\Rt)\Big)\subset M^{p}_{\gamma,0}L^p(0,T),
\end{equation}
with $r<q<+\infty$, which is proven in the Appendix \ref{AppendixB}.


\section{The local Morrey spaces}\label{Sec:Morrey} 
In this section, for the completeness of the paper,  we summarize  some previous results on the local Morrey spaces $M^{p}_{\gamma} (\Rt)$ and $M^{p}_{\gamma,0} (\Rt)$ given in (\ref{def-local-morrey})  and (\ref{cond-dec1}) respectively,  and its time-space version $M^{p}_{\gamma}L^p(0,T)$ and  $M^{p}_{\gamma,0}L^p(0,T)$ defined in (\ref{def-time-space-local-Morrey})  and (\ref{cond-dec-2}) respectively. 

\medskip 

These kind of  local Morrey spaces are strongly lied with the weighted Lebesgue spaces $L^{p}_{w_\gamma}(\Rt)$ which are defined as follows: for $\gamma\geq 0$ we consider the weight 
\begin{equation}\label{peso}
w_\gamma (x)= \frac{1}{(1+\vert x \vert)^\gamma}
\end{equation}
and then for $1< p < +\infty$ we define the space $\ds{L^{p}_{w_\gamma}(\Rt)= L^{p}(w_\gamma\, dx)}$.  Thus, we have the following useful result.
\begin{Lemme}[Lemma $2.1$ of \cite{PFOJ}]\label{Lema-Tech1} Let $0 \leq \gamma < \delta$ and $1<p<+\infty$. 
\begin{enumerate}
\item[$1)$] We have the continuous embedding: $\ds{L^{p}_{w_{\gamma}}(\Rt) \subset M^{p}_{\gamma,0}(\Rt) \subset M^{p}_{\gamma}(\Rt) \subset L^{p}_{w_{\delta}}(\Rt)}$.
\item[$2)$] Moreover, for $0<T<+\infty$ we  also have  the continuous embedding: 
$$L^p\left([0,T], L^{p}_{w_{\gamma}}(\Rt)\right) \subset M^{p}_{\gamma,0}L^p(0,T) \subset M^{p}_{\gamma}L^p(0,T) \subset L^p\left([0,T], L^{p}_{w_{\delta}}(\Rt)\right).$$
\end{enumerate}		
\end{Lemme}	
Thereafter, a second useful result is the following one.
\begin{Lemme}[Lemma $2.1$ of \cite{PFOJ1} and Corollary $2.1$ of \cite{PFOJ}]\label{Lema-Tech2} Let $0 < \gamma < 3$ and $1<p<+\infty$.
\begin{enumerate}
\item[$1)$] The  Riesz transform $\ds{\mathcal{R}_i = \frac{\partial_i}{\sqrt{-\Delta}}}$ is bounded 	on $M^{p}_{\gamma}(\Rt)$ and we have $\ds{\Vert \mathcal{R}_i f \Vert_{M^{p}_{\gamma}} \leq c_{p,\gamma} \Vert f \Vert_{M^{p}_{\gamma}}}$. 
\item[$2)$] The  Hardy-Littlewood maximal function operator $\mathcal{M}$ is also  bounded 	on the space  $M^{p}_{\gamma}(\Rt)$ and we have $\ds{\Vert \mathcal{M}_f \Vert_{M^{p}_{\gamma}} \leq c_{p,\gamma} \Vert f \Vert_{M^{p}_{\gamma}}}$.
\item[$3)$] The points $1)$ and $2)$ also hold for the weighted Lebesgue spaces $\ds{L^{p}_{w_\gamma}(\Rt)}$.
\end{enumerate}		
\end{Lemme}

\section{Characterization of the pressure term}\label{Sec:Presion}
\subsection*{Proof of Proposition \ref{Prop-presion}}

 First, we define $q$ given by the expression 
\begin{equation}\label{Presion}
 q = \sum_{i,j=1}^{3} \mathcal{R}_i \mathcal{R}_j (u_i\, u_j)+ \sum_{i,j,k=1}^{3} \mathcal{R}_i \mathcal{R}_j\left( \partial_i v_k\, \partial_j v_k\right), 
\end{equation}	where, for $9/4<\delta <3$ we have  $\ds{ q \in L^{p/2}([0,T], L^{p/2}_{w_{\delta}}(\Rt))}$. Indeed,  since we have   assumed $\vu \in M^{p}_{\gamma}L^{p}(0,T)$ and $\nabla \otimes \v \in M^{p}_{\gamma}L^{p}(0,T)$, with $0<\gamma< 3/2$,  
then by point $2)$ of Lemma \ref{Lema-Tech1} we get $\ds{\vu \in L^{p}([0,T], L^{p}_{w_{\delta}}(\Rt))}$ and $\ds{\nabla \otimes \v \in L^{p}([0,T], L^{p}_{w_{\delta}}(\Rt))}$. With this information  we are able to  write 
$\ds{\vu \otimes \vu \in L^{p/2}([0,T], L^{p/2}_{w_{\delta}}(\Rt))}$ and $\ds{ \nabla \otimes \v \odot \nabla \otimes \v \in L^{p/2}([0,T], L^{p/2}_{w_{\delta}}(\Rt))}$, and  moreover,  as by point $3)$ of Lemma \ref{Lema-Tech2} the operator $\mathcal{R}_{i}\mathcal{R}_{i}$ is bounded in $\ds{L^{p/2}([0,T], L^{p/2}_{w_{\delta}}(\Rt))}$, then we obtain $\ds{ q \in L^{p/2}([0,T], L^{p/2}_{w_{\delta}}(\Rt))}$. 

\medskip 

Now,  we will prove the identity  $\nabla p= \nabla q$.  For this, let $\varepsilon >0$ (small enough)  and let $\alpha \in \mathcal{C}^{\infty}_{0}(\R)$ be a function such that $\alpha(t)=0$ for $\vert t \vert > \varepsilon$. Moreover, let $\varphi \in \mathcal{C}^{\infty}_{0}(\Rt)$. We may observe that we have $\ds{(\alpha \varphi) \ast \nabla p \in  \mathcal{D}^{'}(]\varepsilon,T-\varepsilon[\times \Rt) }$ and $\ds{(\alpha \varphi) \ast \nabla q \in  \mathcal{D}^{'}(]\varepsilon,T-\varepsilon[\times \Rt) }$ and then, for $t \in ]\varepsilon, T-\varepsilon[$ fix,  we define the expression  $\ds{A_\varepsilon(t)=(\alpha \varphi) \ast \nabla p(t,\cdot)- (\alpha \varphi) \ast \nabla q(t,\cdot) \in \mathcal{D}^{'}(\Rt)}$, where we must verify that we have $A_\varepsilon(t) \in \mathcal{S}^{'}(\Rt)$. We write $\ds{A_\varepsilon(t)=(\alpha \varphi) \ast \nabla p(t,\cdot)- (\alpha \nabla \varphi) \ast  q(t,\cdot)}$. Moreover, since $(\vu, p,  v)$ verify the coupled system (\ref{EickLes-non-stat}) then we have 
$$ \nabla p= -\partial_t \vu +\Delta \vu - div (\vu \otimes \vu) - div (\nabla \otimes \v \odot  \nabla \otimes \v), $$ and thus we obtain 
\begin{equation}\label{A}
\begin{split}
A_\varepsilon(t)=&\left[ (-(\partial_t \alpha) \varphi + \alpha \Delta \varphi )\ast \vu \right](t,\cdot)- \left[(\alpha\ast \nabla \varphi) \ast  (\vu \otimes \vu)\right](t,\cdot)\\
& - \left[(\alpha\ast \nabla \varphi) \ast  (\nabla \otimes \v \odot \nabla \otimes \v )\right](t,\cdot)-\left[(\alpha \nabla \varphi) \ast  q\right](t,\cdot).
\end{split}
\end{equation}
In this identity,  we will prove that each term in the right side  belong to the space $\ds{L^{p/2}_{w_{\delta}} (\Rt)}$ (where $9/4<\delta < 3$). For the first term to the right in (\ref{A}), recall that we have $\vu  \in L^{p}([0,T], L^{p}_{w_{\delta}} (\Rt))$. Moreover, since for $\varphi \in \mathcal{C}^{\infty}_{0}(\Rt)$  and a function $f$ we have the pointwise estimate  $\vert (\varphi \ast f) (x) \vert \leq c_\varphi \,  \mathcal{M}_{f}(x)$ (where $\mathcal{M}$ always denote the Hardy-Littlewood  maximal function operator) then by point $3)$ of Lemma \ref{Lema-Tech2} we obtain that convolution with test functions is a bounded operator on $L^{p}_{w_{\delta}} (\Rt)$. Thus, we have $\ds{\left[(-(\partial_t \alpha) \varphi + \alpha \Delta \varphi )\ast \vu\right](t,\cdot) \in L^{p}_{w_{\delta}} (\Rt)}$. On the other hand, for $9/4<\delta <3$ we have the continuous embedding $\ds{L^{p}_{w_{\delta}} (\Rt) \subset L^{p/2}_{w_{\delta}} (\Rt)}$. Indeed, by definition of the weight $w_\delta(x)$ given by (\ref{peso}) and  using the Cauchy-Schwarz inequalities  we write 
$$ \int_{\Rt} \vert f \vert^{p/2} w_\delta dx = \int_{\Rt} \vert f \vert^{p/2} w_{3/4} w_{\delta-3/4} \leq \left(\int_{\Rt} \vert f \vert^p w_{3/2} dx \right)^{1/2} \left( \int_{\Rt} w_{2\delta - 3/2} dx \right)^{1/2},$$
where, as we have  $9/4 < \delta <3$ then the last integral in the right side convergences. Thus we obtain $\ds{\left[(-(\partial_t \alpha) \varphi + \alpha \Delta \varphi )\ast \vu\right](t,\cdot) \in L^{p/2}_{w_{\delta}} (\Rt)}$.  

\medskip 

For the second and third  terms to the right in (\ref{A}), recall that we have  $\ds{\vu \otimes \vu \in L^{p/2}([0,T], L^{p/2}_{w_{\delta }}(\Rt))}$ and $\ds{ \nabla \otimes \v \odot \nabla \otimes \v \in L^{p/2}([0,T], L^{p/2}_{w_{\delta}}(\Rt))}$, hence, always by the fact that convolution with test functions is a bounded operator on the space $L^{p/2}_{w_{\delta}} (\Rt)$,  we obtain $\ds{\left[(\alpha\ast \nabla \varphi) \ast  (\vu \otimes \vu)\right](t,\cdot) \in L^{p/2}_{w_{\delta}}(\Rt)}$ and $\ds{\left[(\alpha\ast \nabla \varphi) \ast  (\nabla \otimes \v \odot \nabla \otimes \v )\right](t,\cdot) \in L^{p/2}_{w_{\delta}}(\Rt)}$ respectively. 

\medskip 

Finally, for the fourth term to the right in identity (\ref{A}), as we have $\ds{ q \in L^{p/2}([0,T], L^{p/2}_{w_{\delta}}(\Rt))}$ then  we obtain $\ds{\left[(\alpha \nabla \varphi) \ast  q \right](t,\cdot) \in L^{p/2}_{w_{\delta}}(\Rt)}$.

\medskip 

Getting back to the identity (\ref{A}) we get 
$\ds{A_\varepsilon(t) \in L^{p/2}_{w_{\delta}}(\Rt)}$ and then we have  $\ds{A_\varepsilon(t)\in \mathcal{S}^{'}(\Rt)}$.  On the other hand,  since we have $div(\vu)=0$, taking the divergence operator in the first equation of (\ref{EickLes-non-stat}) we obtain $\Delta (p-q)=0$. Then we have  $\ds{\Delta A_\varepsilon(t) =0}$ and since  $\ds{A_\varepsilon(t)\in \mathcal{S}^{'}(\Rt)}$ we get that  $A_\varepsilon(t)$ is a polynomial. But, recalling that  $\ds{A_\varepsilon(t) \in L^{p/2}_{w_{\delta}}(\Rt)}$, we necessary have  $A_\varepsilon(t)=0$. Finally, we use the   approximation of the identity $\ds{\frac{1}{\varepsilon^4}\alpha\left( \frac{t}{\varepsilon}\right) \varphi \left( \frac{x}{\varepsilon}\right)}$  to write $\ds{\nabla(p-q)(t,\cdot)= \lim_{\varepsilon \to 0} A_\varepsilon(t)=0}$.  \finpv

\section{The stationary case}\label{Sec:Stat} 

\subsection{Proof of Theorem  \ref{Th1}} 
All the results stated in this theorem deeply base on the following local estimate, also know as a Cacciopoli-type estimate. 
\begin{Proposition}\label{Prop-Cacciopoli} Let $(\U, p, \V)$ be a  smooth solution of the coupled system (\ref{EickLes}). Let $3\leq p <+\infty$. If $(\U, p) \in L^{p}_{loc}(\Rt)$ and moreover, if $\nabla \otimes \V \in L^{p}_{loc}(\Rt)$, then there exists a constant $c>0$ such that  for all $R \geq 1$ we have: 
\begin{equation}\label{Cacciopoli} 
\begin{split}
&\int_{B_{R/2}} \vert \nabla \otimes \U \vert^2 dx \leq c \left[ \left( \int_{\mathcal{C}(R/2, R)} \vert \U \vert^p dx\right)^{2/p} + \left( \int_{\mathcal{C}(R/2, R)} \vert \nabla \otimes \V \vert^p dx\right)^{2/p} \right. \\
& \left. + \left( \int_{\mathcal{C}(R/2, R)} \vert p \vert^{p/2} dx\right)^{2/p}\right] \times R^{2-9/p} \left( \int_{\mathcal{C}(R/2, R)} \vert \U \vert^p dx \right)^{1/p} + \frac{c}{R^2}\int_{\mathcal{C}(R/2,R)}\vert \U \vert^2 dx.   
\end{split}
\end{equation}
\end{Proposition}	 
\pv  We start by introducing the following cut-off function. Let $\theta \in \mathcal{C}^{\infty}_{0}(\Rt)$ be a positive and radial  function such that $\theta (x)=1$ for $\vert x \vert < 1/2$ and $\theta (x)=0$ for $\vert x \vert \geq 1 $. Then, for $R\geq 1$ we define the function 
\begin{equation}\label{funTheta}
  \theta_R (x)= \theta (x / R).  
\end{equation}
 Remark that this function verifies the following properties: we have $\theta_R (x)=1$ for $\vert x \vert < R /2$, $\theta_R(x)=0$ for $\vert x \vert > R$, and moreover we have $\Vert \nabla \theta_R \Vert_{L^{\infty}} \leq \frac{c}{R} $ and $\Vert \Delta \theta_R \Vert_{L^{\infty}} \leq \frac{c}{R^2}$.  
 
 \medskip 

For $R \geq 1$, we multiply the first equation of the system (\ref{EickLes}) by $\theta_R \U$ and integrating on the ball $B_R$ (since  we have  $supp (\theta_R) \subset B_R$) we obtain: 

\begin{equation}\label{eq01}
\begin{split}
-\int_{B_R} \Delta \U \cdot \theta_R \U dx  + \int_{B_R} \text{div}(\U \otimes \U ) \cdot \theta_R \U dx + \int_{B_R}  \text{div} (\nabla\otimes \V \odot \nabla \otimes \V) \cdot \theta_R \U dx \\
+\int_{B_R}\nabla p \cdot \theta_R \U dx =0. 
\end{split}
\end{equation}
Moreover,  we multiply the second equation of the system (\ref{EickLes}) by $-\theta_R \Delta \V$, then we integrate on the ball $B_R$ to get:  
\begin{equation}\label{eq02}
\int_{B_R}\Delta \V\cdot  \theta_R \Delta \V dx -\int_{B_R} div(\V \otimes \U )\cdot \theta_R \Delta \V dx + \int_{B_R} \vert  \nabla\otimes \V  \vert^2 \V \cdot \theta_R \Delta \V dx =0. 
\end{equation}
At this point remark that as $\U, p$ and $\v$ are smooth functions then all the terms in equations (\ref{eq01}) and (\ref{eq02}) are well-defined. 

\medskip 

Now, we need to study each term in these equations.   We start by equation (\ref{eq01}). For the first term in the left-hand side, by integration by parts we have 
\begin{eqnarray*}
& & 	-\int_{B_R} \Delta \U \cdot \theta_R \U dx =    -\sum_{i,j=1}^{3} \int_{B_R} (\partial^{2}_{j} u_i) ( \theta_R u_i) d =  \sum_{i,j=1}^{3} \int_{B_R} \partial_{j}u_i \partial_{j}(\theta_R u_i) dx \\
&= & \sum_{i,j=1}^{3} \int_{B_R} (\partial_j u_i) (\partial_j \theta_R)    u_i dx   + \sum_{i,j=1}^{3} \int_{B_R} (\partial_j u_i)  \theta_R (  \partial_j u_i) dx \\
& =&  \frac{1}{2} \sum_{i,j=1}^{3} \int_{B_R}  (\partial_j \theta_R) \partial_j (u^{2}_{i})  dx + \int_{B_R} \vert \nabla\otimes \U \vert^2 \theta_R dx\\
	&=& - \frac{1}{2} \int_{B_R} \vert \U \vert^2 \Delta \theta_R dx +   \int_{B_R} \vert \nabla \otimes \U \vert^2 \theta_R dx.  
\end{eqnarray*}
For the second term in the left-hand side of (\ref{eq01}), by integration by parts and moreover, as we have $div(\U)=0$,  we can write 
\begin{eqnarray*}
& & \int_{B_R} div(\U \otimes \U ) \cdot \theta_R \U dx = \sum_{i,j=1}^{3} \int_{B_R} \partial_j (u_i u_j) \theta_R u_i dx  \\
&=& -\sum_{i,j=1}^{3} \int_{B_R} u_i u_j (\partial_{j}\theta_R)  u_i dx- \sum_{i,j=1}^{3} \int_{B_R} u_i u_j \theta_R  (\partial_j   u_i) dx\\
&=& - \int_{B_R} \vert \U \vert^2 (\U \cdot \nabla \theta_R) dx - \frac{1}{2} \sum_{i,j=1}^{3} \int_{B_R} u_j \theta_R \partial_{j} (u^{2}_{i}) dx\\
&=& - \int_{B_R} \vert \U \vert^2 (\U \cdot \nabla \theta_R) dx + \frac{1}{2} \sum_{i,j=1}^{3} \int_{B_R} \partial_{j} (u_j \theta_R) u^{2}_{i} dx \\
&=& - \int_{B_R} \vert \U \vert^2 (\U \cdot \nabla \theta_R) dx  + \frac{1}{2} \int_{B_R} (\U \cdot \nabla \theta_R) \vert \vu\vert^2  dx \\
&=& - \frac{1}{2} \int_{B_R} \vert \vu \vert^2 (\U \cdot \nabla\theta_R) dx.  
\end{eqnarray*}
In order to study the third term in the left-hand side of  (\ref{eq01}), we need the following technical identity: 
$$  div( \nabla\otimes \V \odot \nabla\otimes \V  )= \nabla \left( \frac{1}{2} \vert \nabla \otimes \V \vert^2 \right)+ \Delta \V  ( \nabla \otimes \V) .$$ 
Indeed, recall that   for $i=1,2,3$ each component of the vector field   $div( \nabla \otimes \V \odot \nabla\otimes \V  )$ is given by 
\begin{eqnarray*}
& & (div( \nabla\otimes \V \odot \nabla \otimes \V  ))_i= \sum_{j,k=1}^{3} \partial_j (\partial_i v_k \, \partial_j v_k)= \sum_{j,k=1}^{3} \partial_j (\partial_i v_k) \partial_j v_k + \sum_{j,k=1}^{3} \partial_{i} v_k \partial^{2}_j v_k\\
&=& \sum_{j,k=1}^{3} \partial_i (\partial_j v_k)  \partial_j v_k + \sum_{k=1}^{3} \partial_i v_k \, \Delta v_k= \partial_i \left( \frac{1}{2} \sum_{j,k=1}^{3} (\partial_j v_k)^2\right)+\sum_{k=1}^{3} \Delta v_k\,  \partial_i v_k\\
&=& \partial_i \left( \frac{1}{2} \vert \nabla \otimes \V \vert^2\right)+ (\Delta \V (\nabla \otimes \V))_{i}. 
\end{eqnarray*}

With this identity at hand, we get back to the third term in the left-hand side in (\ref{eq01}) and,   by integration by parts and the fact that $div(\U)=0$,    we  write 
\begin{eqnarray*}
& & 	\int_{B_R}  div (\nabla \otimes \V \odot \nabla \otimes \V) \cdot \theta_R \U dx = \sum_{i=1}^{3} \int_{B_R} \partial_i \left( \frac{1}{2} \vert \nabla \otimes \V \vert^2 \right) \theta_R U_i dx \\
& & + \sum_{i,j=1}^{3}\int_{B_R} \Delta V_j (\partial_i V_j)  \theta_R u_i dx  = - \frac{1}{2} \int_{B_R} \vert \nabla \otimes \V \vert^2 (\U \cdot \nabla \theta_R)dx \\
& & + \sum_{i,j=1}^{3}\int_{B_R} \Delta v_j (\partial_i v_j)  \theta_R u_i dx.  
\end{eqnarray*}
Finally, for the fourth term in the left-hand side in (\ref{eq01}), always by integration by parts and since  $div(\U)=0$ we have 
\begin{equation*}
\int_{B_R}\nabla p \cdot \theta_R \U dx= \sum_{i=1}^{3} \int_{B_R} (\partial_i p ) \theta_R u_i dx = - \int_{B_R} p (\U \cdot \nabla \theta_R) dx . 
\end{equation*}
Once we dispose of these identities, we  get back to equation (\ref{eq01})  and then we obtain 
\begin{equation}\label{eq03}
\begin{split}
\int_{B_R} \vert \nabla \otimes \U \vert^2 \theta_R dx = \int_{B_R} \left( \frac{\vert \U \vert^2}{2} + \frac{\vert \nabla \otimes \V \vert^2}{2} + p \right)(\U \cdot \nabla \theta_R)dx \\
+ \frac{1}{2} \int_{B_R} \vert \U \vert^2 \Delta \theta_R dx 
-  \sum_{i,j=1}^{3}\int_{B_R} \Delta v_j (\partial_i v_j)  \theta_R u_i dx.
\end{split}
\end{equation}
We study now the terms in the left-hand side in equation (\ref{eq02}). For the first term we write  directly
$$ \int_{B_R} \Delta \V \cdot \theta_R \Delta \V dx = \int_{B_R} \vert \Delta \V \vert^2 \theta_R dx.  $$
For the second term, integrating by parts and as $div(\U)=0$  then we get
\begin{eqnarray*}
& &- \int_{B_R} div(\V \otimes \U) \cdot \theta_R \Delta \V dx =  - \sum_{i,j=1}^{3} \int_{B_R} \partial_{j}(v_i u_j) \theta_R  \Delta v_i \\
&=-& \sum_{i,j=1}^{3} \int_{B_R} (\partial_j v_i) u_j \theta_R \Delta v_i =  - \sum_{i,j=1}^{3} \int_{B_R} \Delta v_i (\partial_j v_i) \theta_R u_j dx. 
\end{eqnarray*}
For the third term we write 
\begin{equation*}
\int_{B_R} \vert \nabla \otimes \V \vert^2 \V \cdot \theta_R \Delta \V dx = \sum_{i=1}^{3} \int_{B_R} \vert \nabla \otimes \V \vert^2 v_i \theta_R \Delta v_i dx = \int_{B_R} \vert \nabla \otimes \V \vert^2 (\V \cdot \Delta \V) \theta_R dx. 
\end{equation*}
Thus, with these identities at hand, from equation (\ref{eq02}) we obtain: 
\begin{equation}\label{eq04}
\int_{B_R} \vert \Delta \V \vert^2 \theta_R dx = \sum_{i,j=1}^{3} \int_{B_R} \Delta v_i (\partial_j v_i) \theta_R u_j dx - \int_{B_R} \vert \nabla \otimes \V \vert^2 (\V \cdot \Delta \V) \theta_R dx. 
\end{equation}
Now, adding the equations (\ref{eq03}) and (\ref{eq04}) we get 
\begin{eqnarray*}
& & \int_{B_R} \vert \nabla\otimes \U \vert^2 \theta_R dx + \int_{B_R} \vert \Delta \V \vert^2 \theta_R dx = 	\int_{B_R} \left( \frac{\vert \U \vert^2}{2}+ \frac{\vert \nabla \otimes \V \vert^2}{2}+ p \right)(\U \cdot \nabla \theta_R)dx \\
& & +  \int_{B_R} \frac{\vert \U \vert^2 }{2} \Delta \theta_R dx  \underbrace{ -  \sum_{i,j=1}^{3}\int_{B_R} \Delta v_j (\partial_i v_j)  \theta_R u_i dx + \sum_{i,j=1}^{3} \int_{B_R} \Delta v_i (\partial_j v_i) \theta_R u_j dx}_{(a)}\\
& & -  \int_{B_R} \vert \nabla \otimes \V \vert^2 (\V \cdot \Delta \V) \theta_R dx, 
\end{eqnarray*}
but, we may observe that we have $(a)=0$ and then we write 
\begin{eqnarray*}
	& & \int_{B_R} \vert \nabla\otimes \U \vert^2 \theta_R dx + \int_{B_R} \vert \Delta \V \vert^2 \theta_R dx = 	\int_{B_R} \left( \frac{\vert \U \vert^2}{2} + \frac{\vert \nabla \otimes \V \vert^2}{2} + p \right)(\U \cdot \nabla \theta_R)dx \\
	& &  +  \int_{B_R} \frac{\vert \U \vert^2 }{2} \Delta \theta_R dx -  \int_{B_R} \vert \nabla \otimes \V \vert^2 (\V \cdot \Delta \V) \theta_R dx. 
\end{eqnarray*}
Moreover, the last term is  estimated as follows:  
$$  -  \int_{B_R} \vert \nabla \otimes \V \vert^2 (\V \cdot \Delta \V) \theta_R dx \leq \int_{B_R} \vert \Delta \V \vert^2 \theta_R dx. $$  
Indeed, recall that by hypothesis  we have  $\vert \V \vert^{2}=1$ and then  we get  $ \frac{1}{2} \Delta \vert \V \vert^2 =0$. Thus, we can write  

\begin{equation}\label{Iden}	
\begin{split}
&  -\vert \nabla\otimes \V \vert^2 = - \sum_{i,j=1}^{3} (\partial_i v_j)^2 = - \sum_{i,j=1}^{3} (\partial_i v_i)^2 = - \sum_{i,j=1}^{3} (\partial_i v_j)^2+ \frac{1}{2} \Delta \vert \V \vert^2 \\
= & - \sum_{i,j=1}^{3} (\partial_i v_j)^2  + \frac{1}{2} \sum_{i,j=1}^{2} \partial^{2}_{i} (v^{2}_{i})  - \sum_{i,j=1}^{3} (\partial_i v_j)^2  + \sum_{i,j=1}^{3} \partial_j \left( \frac{1}{2} \partial_i  v^{2}_{i}\right)\\
= & - \sum_{i,j=1}^{3} (\partial_i v_j)^2  + \sum_{i,j=1}^{3} \partial_j (v_i \partial_j v_i )  - \sum_{i,j=1}^{3} (\partial_i v_j)^2 + \sum_{i,j=1}^{3} (\partial_j v_i)^2 + \sum_{i,j=1}^{3} v_i \partial^{2}_{j} v_i\\
=&  \sum_{i,j=1}^{3} v_i \partial^{2}_{j} v_i =   \V \cdot \Delta \V.
\end{split} 
\end{equation}
With the identity $\ds{-\vert \nabla \otimes \V \vert^2=\V \cdot \Delta \V }$ at hand,  and moreover, as we have $\theta_R \geq 0$ and  as we have $\vert \v \vert^2=1$,  we obtain 
\begin{equation*}
  -  \int_{B_R} \vert \nabla \otimes \V \vert^2 (\V \cdot \Delta \V) \theta_R dx = \int_{B_R} \vert \V \cdot \Delta \V\vert^2 \theta_R dx \leq \int_{B_R} \vert \V \vert^2 \vert \Delta \V \vert^2 \theta_R dx \leq \int_{B_R} \vert \Delta \V \vert^2 \theta_R dx.   
\end{equation*}
Once we have  this estimate then we can write 
\begin{eqnarray*} 
	& & \int_{B_R} \vert \nabla \otimes \U \vert^2 \theta_R dx + \int_{B_R} \vert \Delta \V \vert^2 \theta_R dx  \leq 	\int_{B_R} \left( \frac{\vert \U \vert^2}{2} + \frac{\vert \nabla \otimes \V \vert^2}{2} + p \right)(\U \cdot \nabla \theta_R)dx \\
	& & +  \int_{B_R} \frac{\vert \U \vert^2 }{2} \Delta \theta_R dx +  \int_{B_R} \vert \Delta \V \vert^2 \theta_R dx, 
\end{eqnarray*}
hence we get 
\begin{equation*}
\int_{B_R} \vert \nabla\otimes \U \vert^2 \theta_R dx \leq 	\int_{B_R} \left( \frac{\vert \U \vert^2}{2} + \frac{\vert \nabla \otimes \V \vert^2}{2} + p \right)(\U \cdot \nabla \theta_R)dx+  \int_{B_R} \frac{\vert \U \vert^2 }{2} \Delta \theta_R dx. 
\end{equation*}
Recalling that we have $\theta_R(x)=1$ for $\vert x \vert < R/2$, then we obtain $$\ds{\int_{B_{R/2}} \vert \nabla \otimes \U \vert^2 dx \leq \int_{B_R} \vert \nabla \otimes \U \vert^2 \theta_R dx},$$ and from the previous inequality we are able to write 
\begin{equation*}
\int_{B_{R/2}} \vert \nabla \otimes \U \vert^2 dx \leq \int_{B_R} \left( \frac{\vert \U \vert^2}{2} + \frac{\vert \nabla \otimes \V \vert^2}{2} + p \right)(\U \cdot \nabla \theta_R)dx+  \int_{B_R} \frac{\vert \U \vert^2 }{2} \Delta \theta_R dx.
\end{equation*}
Moreover, recalling that we have $supp (\nabla \theta_R) \subset \mathcal{C}(R/2, R)$ and  $supp (\Delta \theta_R) \subset \mathcal{C}(R/2, R)$, then we obtain the following estimate 
\begin{equation}\label{eq05}
\begin{split}
 \int_{B_{R/2}} \vert \nabla \otimes \U \vert^2 dx & \leq \int_{\mathcal{C}(R/2, R) } \left( \frac{\vert \U \vert^2}{2} + \frac{\vert \nabla \otimes \V \vert^2}{2} + p \right)(\U \cdot \nabla \theta_R)dx+  \int_{\mathcal{C}(R/2, R) } \frac{\vert \U \vert^2 }{2} \Delta \theta_R dx \\
&\leq  \int_{\mathcal{C}(R/2, R) } \frac{\vert \U \vert^2}{2} (\U \cdot \nabla \theta_R) dx + \int_{\mathcal{C}(R/2, R) } \frac{\vert \nabla\otimes \V \vert^2}{2} (\U \cdot \nabla \theta_R) dx \\
& + \int_{\mathcal{C}(R/2, R) } p (\U \cdot \nabla \theta_R) dx +  \int_{\mathcal{C}(R/2, R) } \frac{\vert \U \vert^2 }{2} \Delta \theta_R dx= I_1+I_2+I_3+I_4.
\end{split}    
\end{equation}

From this estimate we will derive the desired inequality (\ref{Cacciopoli}) and for this we will study each term $I_i$ for $i=1, \cdots 4$.  For the term $I_1$, by the H\"older inequalities (with $1=2/p + 1/q$ )  and moreover, as we have $\Vert \nabla \theta_R \Vert_{L^{\infty}} \leq c / R$, we get  
\begin{eqnarray*}
I_1& \leq & \int_{\mathcal{C}(R/2, R)} \vert \U \vert^2 \vert \U \cdot \nabla \theta_R \vert dx \leq \left( \int_{\mathcal{C}(R/2, R)} \vert \U \vert^p dx\right)^{2/p} \left( \int_{\mathcal{C}(R/2, R)} \vert \U \cdot \nabla\theta_R \vert^q  dx\right)^{1/q}  \\
& \leq & \left( \int_{\mathcal{C}(R/2, R)} \vert \U \vert^p dx\right)^{2/p} \, \frac{c}{R} \left( \int_{\mathcal{C}(R/2, R)} \vert \U  \vert^q  dx\right)^{1/q} .
\end{eqnarray*}
But, since we have $3 \leq p < +\infty $ and $1=2/p + 1/q$ then  the parameter $q$ verifies  $ q \leq  3 \leq p$ and thus, for the last expression  we can write 
\begin{equation*}
 \frac{c}{R} \left( \int_{\mathcal{C}(R/2, R)} \vert \U  \vert^q  dx\right)^{1/q} \leq \frac{c}{R} R^{3(1/q - 1/p)} \left( \int_{\mathcal{C}(R/2, R)} \vert \U \vert^p dx \right)^{1/p} 
 \leq c \, R^{2-9/p} \left( \int_{\mathcal{C}(R/2, R)} \vert \U \vert^p dx \right)^{1/p},    
 \end{equation*}
hence we have 
\begin{equation}\label{I1} 
I_1 \leq c \left( \int_{\mathcal{C}(R/2, R)} \vert \U \vert^p dx\right)^{2/p}  \,   R^{2-9/p} \left( \int_{\mathcal{C}(R/2, R)} \vert \U \vert^p dx \right)^{1/p}. 
\end{equation}
Following the same computations, the terms $I_2$  and $I_3$ are estimated as follows: 
\begin{equation}\label{I2}
I_2 \leq c \left( \int_{\mathcal{C}(R/2, R)} \vert \nabla\otimes \V \vert^p dx\right)^{2/p}  \,   R^{2-9/p} \left( \int_{\mathcal{C}(R/2, R)} \vert \U \vert^p dx \right)^{1/p},
\end{equation} and 
\begin{equation}\label{I3}
I_3 \leq c \left( \int_{\mathcal{C}(R/2, R)} \vert p \vert^{p/2} dx\right)^{2/p}  \,   R^{2-9/p} \left( \int_{\mathcal{C}(R/2, R)} \vert \U \vert^p dx \right)^{1/p}. 
\end{equation}
Finally, for the term $I_4$, always by the H\"older inequalities, with $1=2/p+1/q$,  by the fact that $\Vert \Delta \theta_R \Vert_{L^{\infty}} \leq c / R^2$  we obtain 
$$ I_4  \leq  c\int_{\mathcal{C}(R/2, R)} \vert \U \vert^2 \vert \Delta \theta_R \vert dx \leq \frac{c}{R^2} \int_{\mathcal{C}(R/2, R)} \vert \U \vert^2 dx. $$
With the estimates, we get back to the inequality (\ref{eq05}) to obtain the desired estimate (\ref{Cacciopoli}).  Proposition \ref{Prop-Cacciopoli} is verified.\finpv 

We have now all the tools to prove the Theorem \ref{Th1}. Let $(\U, p, \V)$ be a smooth solution of (\ref{EickLes}) such that $\V$ verifies (\ref{cond-derV}) and moreover  $\U\in M^{p}_{\gamma,0}(\Rt)$ and $\nabla\otimes \V \in M^{p}_{\gamma}(\Rt)$. Then, by Proposition \ref{Prop-presion} we have the identity $\nabla p = \nabla q$, where  $q$ is given in formula (\ref{Presion}); and from now on we will consider the equation (\ref{EickLes}) with the term $\nabla q$ instead of the term $\nabla p$. 

\medskip 

As we have $\U\in M^{p}_{\gamma,0}(\Rt)$ and $\nabla\otimes \V \in M^{p}_{\gamma}(\Rt)$, then by definition of the local Morrey spaces $M^{p}_{\gamma,0}(\Rt)$ and  $ M^{p}_{\gamma}(\Rt)$, given in  (\ref{cond-dec1}) and (\ref{def-local-morrey}) respectively,  we obtain $\U \in L^{p}_{loc}(\Rt)$ and $\nabla\otimes \V \in L^{p}_{loc}(\Rt)$.  Thus, we  assume now $3\leq p <+\infty$; and  by Proposition \rd{\ref{Prop-Cacciopoli}} the local estimate (\ref{Cacciopoli}) holds true for all $R\geq 1$.  Then, for $0<\gamma$ we write


\begin{eqnarray*}
\int_{B_{R/2}} \vert \nabla \otimes \U \vert^2 dx & \leq &   \frac{c}{R^2}\int_{\mathcal{C}(R/2,R)}\vert \U \vert^2 dx + \frac{c}{R^{\frac{2}{p}\gamma}} \left[ \left( \int_{\mathcal{C}(R/2, R)} \vert \U \vert^p dx\right)^{2/p} + \left( \int_{\mathcal{C}(R/2, R)} \vert \nabla \otimes \V \vert^p dx\right)^{2/p}\right.\\
& & \left. +  \left( \int_{\mathcal{C}(R/2, R)} \vert q \vert^{p/2} dx\right)^{2/p}\right]\times R^{\frac{2}{p}\gamma+2- \frac{9}{p} } \left( \int_{\mathcal{C}(R/2, R)} \vert \U \vert^p dx \right)^{1/p}\\
&=&  \frac{c}{R^2}\int_{\mathcal{C}(R/2,R)}\vert \U \vert^2 dx+ c\left[ \left(\frac{1}{R^{\gamma}} \int_{\mathcal{C}(R/2, R)} \vert \U \vert^p dx\right)^{2/p} + \left(\frac{1}{R^{\gamma}} \int_{\mathcal{C}(R/2, R)} \vert \nabla \otimes \V \vert^p dx\right)^{2/p}\right.\\
& & \left. +  \left(\frac{1}{R^{\gamma}} \int_{\mathcal{C}(R/2, R)} \vert q  \vert^{p/2} dx\right)^{2/p}\right]\times  R^{\frac{2}{p}\gamma+2-\frac{9}{p}} \left(  \int_{\mathcal{C}(R/2, R)} \vert \U \vert^p dx \right)^{1/p}.\\
&=& \frac{c}{R^2}\int_{\mathcal{C}(R/2,R)}\vert \U \vert^2 dx+ c\left[ \left(\frac{1}{R^{\gamma}} \int_{\mathcal{C}(R/2, R)} \vert \U \vert^p dx\right)^{2/p} + \left(\frac{1}{R^{\gamma}} \int_{\mathcal{C}(R/2, R)} \vert \nabla \otimes \V \vert^p dx\right)^{2/p}\right.\\
& & \left. +  \left(\frac{1}{R^{\gamma}} \int_{\mathcal{C}(R/2, R)} \vert q  \vert^{p/2} dx\right)^{2/p}\right]\times  R^{\frac{3}{p}\gamma+2-\frac{9}{p}} \left( \frac{1}{R^{\gamma}} \int_{\mathcal{C}(R/2, R)} \vert \U \vert^p dx \right)^{1/p}.\\
\end{eqnarray*}
At this point, recalling that by (\ref{eta}) we define $\ds{\eta(\gamma,p)= \frac{\gamma}{p}  - \frac{3}{p} + \frac{2}{3}}$, then we have  $\ds{ \frac{3}{p}\gamma  +2- \frac{9}{p} = 3\eta(\gamma,p)}$, and we obtain 

\begin{equation}\label{estim-aux}
\begin{split}
\int_{B_{R/2}} \vert \nabla \otimes \U \vert^2 dx & \leq  \frac{c}{R^2}\int_{\mathcal{C}(R/2,R)}\vert \U \vert^2 dx+ c\left[\left(\frac{1}{R^{\gamma}} \int_{\mathcal{C}(R/2, R)} \vert \U \vert^p dx\right)^{2/p} + \left(\frac{1}{R^{\gamma}} \int_{\mathcal{C}(R/2, R)} \vert \nabla \otimes \V \vert^p dx\right)^{2/p}\right.\\ 
& \left. +\left(\frac{1}{R^{\gamma}} \int_{\mathcal{C}(R/2, R)} \vert q  \vert^{p/2} dx\right)^{2/p}\right]\times  R^{3\eta(\gamma,p)} \left( \frac{1}{R^{\gamma}} \int_{\mathcal{C}(R/2, R)} \vert \U \vert^p dx \right)^{1/p}.
\end{split}
\end{equation}  
Here, as we have $\U \in M^{p}_{\gamma,0}(\Rt)$ and $\nabla \otimes \V \in M^{p}_{\gamma}(\Rt)$,  for all $R\geq 1$ we have the uniformly bound 
\[ \left(\frac{1}{R^{\gamma}} \int_{\mathcal{C}(R/2, R)} \vert \U \vert^p dx\right)^{2/p} + \left(\frac{1}{R^{\gamma}} \int_{\mathcal{C}(R/2, R)} \vert \nabla \otimes \V \vert^p dx\right)^{2/p} \leq c  \left( \Vert \U \Vert^{2}_{M^{p}_{\gamma}}+   \Vert \nabla\otimes \V \Vert^{2}_{M^{p}_{\gamma}} \right).\]  
Moreover, in order to estimate the expression  $\ds{\left(\frac{1}{R^{\gamma}} \int_{\mathcal{C}(R/2, R)} \vert q  \vert^{p/2} dx\right)^{2/p}}$, we recall that the term $q$ is defined through $\U$ and $\nabla \otimes \V$ in (\ref{Presion}), and then, 
setting the parameter $0<\gamma$ as $0<\gamma<3$, by the  point $1)$ of Lemma \ref{Lema-Tech2} we also can  write  
\begin{equation}\label{estimP}
\left(\frac{1}{R^{\gamma}} \int_{\mathcal{C}(R/2, R)} \vert q  \vert^{p/2} dx\right)^{2/p} \leq  \Vert q \Vert_{M^{p/2}_{\gamma}} \leq c  \left( \Vert \U \Vert^{2}_{M^{p}_{\gamma}}+   \Vert \nabla\otimes \V \Vert^{2}_{M^{p}_{\gamma}} \right).   
\end{equation}

Getting back to (\ref{estim-aux}),  we have the estimate 
\begin{equation*}
\begin{split}
\int_{B_{R/2}} \vert \nabla \otimes \U \vert^2 dx \leq \frac{c}{R^2}\int_{\mathcal{C}(R/2,R)}\vert \U \vert^2 dx+  c \left( \Vert \U \Vert^{2}_{M^{p}_{\gamma}}+   \Vert \nabla\otimes \V \Vert^{2}_{M^{p}_{\gamma}} \right) \, R^{3\eta(\gamma,p)} \left( \frac{1}{R^{\gamma}} \int_{\mathcal{C}(R/2, R)} \vert \U \vert^p dx \right)^{1/p},
\end{split}
\end{equation*}
where, we still must  study the first  term in the right. Precisely,   as $ R \geq 1$  we have 
\begin{equation}\label{lim-1}
 \frac{c}{R^2}\int_{\mathcal{C}(R/2,R)}\vert \U \vert^2 dx  \leq c\, R^{2 \eta(\gamma,p)} \left( \frac{1}{R^\gamma} \int_{\mathcal{C}(R/2,R)} \vert \U \vert^p\right)^{2/p}. 
\end{equation}
Indeed, we write 
\begin{equation*}
\frac{c}{R^2} \int_{\mathcal{C}(R/ 2, R)}\vert \U \vert^{2}dx \leq c R^{6(1/2-1/p)-2} \left( \int_{\mathcal{C}(R/2,R)} \vert \U \vert^p\right)^{2/p} \leq c R^{\,6(1/2-1/p)-2 + 2 \gamma /p} \left(\frac{1}{R^\gamma} \int_{\mathcal{C}(R/2,R)} \vert \U \vert^p\right)^{2/p}, 
\end{equation*}
where, always by (\ref{eta}) we can write $\ds{ 6(1/2-1/p)-2 + 2 \gamma /p \leq 2(\gamma/p-3/p+1/2) \leq 2 \eta(\gamma,p)}$, hence the estimate (\ref{lim-1}) follows. 

\medskip 

Thus, we obtain the following estimate 

\begin{equation}\label{estim-base}
\begin{split}
\int_{B_{R/2}} \vert \nabla \otimes \U \vert^2 dx \leq & \,\, c\, R^{2 \eta(\gamma,p)} \left( \frac{1}{R^\gamma} \int_{\mathcal{C}(R/2,R)} \vert \U \vert^p\right)^{2/p}\\
& +  c \left( \Vert \U \Vert^{2}_{M^{p}_{\gamma}}+   \Vert \nabla\otimes \V \Vert^{2}_{M^{p}_{\gamma}} \right) \, R^{3\eta(\gamma,p)} \left( \frac{1}{R^{\gamma}} \int_{\mathcal{C}(R/2, R)} \vert \U \vert^p dx \right)^{1/p},
\end{split}
\end{equation}
and now, we will consider the cases when $\eta(\gamma,p)\leq 0$ and $\eta(\gamma,p)>0$ separately. 

\begin{enumerate}
\item[$1)$] \textbf{The case when  $\eta(\gamma,p)\leq 0$}. Here, as $R \geq 1$  then we have $\ds{R^{2 \eta(\gamma,p) } \leq 1}$ and $\ds{R^{3 \eta(\gamma,p) } \leq 1}$. Thus,  by the estimate (\ref{estim-base}) we can write 
\begin{equation*}
\int_{B_{R/2}} \vert \nabla \otimes \U \vert^2 dx \leq c \left( \frac{1}{R^\gamma} \int_{\mathcal{C}(R/2,R)} \vert \U \vert^p\right)^{2/p}  +  c \left( \Vert \U \Vert^{2}_{M^{p}_{\gamma}}+   \Vert \nabla\otimes \V \Vert^{2}_{M^{p}_{\gamma}} \right) \,  \left( \frac{1}{R^{\gamma}} \int_{\mathcal{C}(R/2, R)} \vert \U \vert^p dx \right)^{1/p}.
\end{equation*}
Moreover, as  $\U \in M^{p}_{\gamma,0}(\Rt)$, taking the limit when $R\to +\infty$ in each side of the estimate above we obtain $\ds{\int_{\Rt} \vert \nabla\otimes \U \vert ^2 dx =0}$ and thus $\U$ is a constant vector. But, always by the information  $\U \in M^{p}_{\gamma,0}(\Rt)$  we necessary have the identity $\U=0$. \\
\item[$2)$] \textbf{The case when  $\eta(\gamma,p)>0$.}  Here, always as $R \geq 1$ then we have $\ds{R^{2 \eta(\gamma,p)}\leq R^{6\eta(\gamma,p)}}$; and thus,   by the estimate (\ref{estim-base}) we  write  now 
\begin{equation*}
\begin{split}
\int_{B_{R/2}} \vert \nabla \otimes \U \vert^2 dx \leq & \,\, c\, R^{6 \eta(\gamma,p)} \left( \frac{1}{R^\gamma} \int_{\mathcal{C}(R/2,R)} \vert \U \vert^p\right)^{2/p}\\
&  +  c \left( \Vert \U \Vert^{2}_{M^{p}_{\gamma}}+   \Vert \nabla\otimes \V \Vert^{2}_{M^{p}_{\gamma}} \right) \, R^{3\eta(\gamma,p)} \left( \frac{1}{R^{\gamma}} \int_{\mathcal{C}(R/2, R)} \vert \U \vert^p dx \right)^{1/p}\\
\leq & c\,\left[ R^{3 \eta(\gamma,p)} \left( \frac{1}{R^\gamma} \int_{\mathcal{C}(R/2,R)} \vert \U \vert^p\right)^{1/p} \right]^{2} \\
&+  c \left( \Vert \U \Vert^{2}_{M^{p}_{\gamma}}+   \Vert \nabla\otimes \V \Vert^{2}_{M^{p}_{\gamma}} \right) \, \left[ R^{3\eta(\gamma,p)} \left( \frac{1}{R^{\gamma}} \int_{\mathcal{C}(R/2, R)} \vert \U \vert^p dx \right)^{1/p} \right]. 
\end{split}
\end{equation*}
Hence, as $\U \in \U \in M^{p}_{\gamma,0}(\Rt)$, and moreover, assuming the supplementary  decaying condition  (\ref{cond-dec}),  have the identity $\U=0$.  
\end{enumerate}

Until now we have proven that $\U=0$ and then it remains to prove the identities $\nabla \otimes \V =0$ and $q=0$. We start by proving that $\nabla \otimes \V =0$. As $\U=0$ then by  (\ref{EickLes}) we have that $\V$ solves the following elliptic equation 
\begin{equation*}
-\Delta \V - \vert \nabla\otimes\V\vert^2 \V=0.     
\end{equation*}
In this equation, we multiply by $\theta_R ((x \cdot \nabla) \V )$, where for $R\geq 1$ the cut-off function $\theta_R(x)$ was defined in (\ref{funTheta}), and integrating on the ball $B_R$  by \bl{\cite{HaLiZ}}, page $6$, we have the following local estimate:
\begin{equation}\label{estimV}
\int_{B_{R/2}} \vert \nabla\otimes \V \vert^2 dx \leq c \int_{\mathcal{C}(R/2, R)} \vert  \nabla\otimes \V \vert^2 dx.    
\end{equation}
 Now, recall that  $\V$ verifies (\ref{cond-derV}) and then we have 
$$ \int_{B_{R/2}} \vert \nabla\otimes \V \vert^2 dx \leq c\, \sup_{R\geq 1} \int_{\mathcal{C}(R/2, R)} \vert  \nabla\otimes \V \vert^2 dx<+\infty,$$ hece we obtain  $\ds{\int_{\Rt}\vert \nabla\otimes \V \vert^2 dx <+\infty }$. With this information, we get back to (\ref{estimV}) and taking the limit when $R\to +\infty$ we get $\ds{\int_{\Rt}\vert \nabla\otimes \V \vert^2 dx =0}$. Hence we have $\nabla \otimes \V =0$. Once we have the identities $\U=0$ and $\nabla \otimes \V =0$, the identity $q=0$ follows directly from the estimate (\ref{estimP}). Finally, always by the identity $\nabla p = \nabla q$ given by Proposition \ref{Prop-presion}, we conclude that $p$ is a constant vector. Theorem  \rd{\ref{Th1}} is proven. \finpv

\subsection{Proof of Corollary \ref{Corollary-NS}}

We observe first that if  $(\U, p)$  is a smooth solution of the equations  (\ref{NS}),  then, for a constant vector field $\V \in \mathbb{S}^{n-1}$ the triplet $(\U,p,\V)$ is also a smooth solution  of the coupled system (\ref{EickLes}). Thus, assuming the velocity $\U$ verifies $\U \in M^{p}_{\gamma,0}(\Rt)$,  where $0<\gamma<3 \leq p <+\infty$ are such that $\eta(\gamma,p)\leq 0$,  and moreover, as we have $\nabla \otimes \V=0$ and consequently $\vec{\nabla} \otimes \V \in M^{p}_{\gamma}(\Rt)$, the result stated in this corollary directly follows from Theorem \ref{Th1}. \finpv 

\subsection{Proof of Proposition \ref{Prop:NS}}  
As mentioned above, the stationary Navier-Stokes equations (\ref{NS}) can be observed as a particular of the coupled Ericksen-Leslie system (\ref{EickLes}) when the unitary vector field $\V$ is a constant vector. Then,
the Proposition \ref{Prop-presion}  holds true for the equations (\ref{NS}), and we write the term $\nabla q$ instead of the term $\nabla p$, where, $\nabla \otimes \V=0$, the term  $q$ is given by the identity $\ds{q=\sum_{i,j=1}^{3} \mathcal{R}_i \mathcal{R}_j (u_i\, u_j)}$. 

\medskip

We also observe that  Proposition \ref{Prop-Cacciopoli}  holds true for the equations (\ref{NS}), and, always  as we have $\nabla \otimes \V=0$, then  we are able to write  the following estimate:
\begin{equation*}
\begin{split}
\int_{B_{R/2}} \vert \nabla \otimes \U \vert^2 dx \leq & \, c \left[ \left( \int_{\mathcal{C}(R/2, R)} \vert \U \vert^p dx\right)^{2/p}   + \left( \int_{\mathcal{C}(R/2, R)} \vert q \vert^{p/2} dx\right)^{2/p}\right] \\
&\times R^{2-9/p} \left( \int_{\mathcal{C}(R/2, R)} \vert \U \vert^p dx \right)^{1/p} + \frac{c}{R^2}\int_{\mathcal{C}(R/2,R)}\vert \U \vert^2 dx.   
\end{split}
\end{equation*}
Hence, following the same computations performed in the estimate (\ref{estim-base}) we obtain
\begin{equation*}
\begin{split}
\int_{B_{R/2}} \vert \nabla \otimes \U \vert^2 dx \leq & \,\, c\, R^{2 \eta(\gamma,p)} \left( \frac{1}{R^\gamma} \int_{\mathcal{C}(R/2,R)} \vert \U \vert^p\, dx\right)^{2/p}\\
& +  c \left( \Vert \U \Vert^{2}_{M^{p}_{\gamma}} \right) \, R^{3\eta(\gamma,p)} \left( \frac{1}{R^{\gamma}} \int_{\mathcal{C}(R/2, R)} \vert \U \vert^p dx \right)^{1/p},
\end{split}
\end{equation*}
and moreover, recalling the definition of the quantity $\Vert \U \Vert_{M^{p}_{\gamma}}$ given in (\ref{def-local-morrey}) we finally have the following estimate:

\begin{equation}\label{estim-base-NS}
\int_{B_{R/2}} \vert \nabla \otimes \U \vert^2 dx \leq  \,\, c\, R^{2 \eta(\gamma,p)} \Vert \U \Vert^{2}_{M^{p}_{\gamma}}+ c\, R^{3\eta(\gamma,p)} \Vert \U \Vert^{3}_{M^{p}_{\gamma}}.
\end{equation}
In this estimate, we will distinguish two cases when $\eta(\gamma,p)<0$ and when $\eta(\gamma,p)=0$.
\begin{enumerate}
\item[$1)$] The case  	$\eta(\gamma,p)<0$. Here, in each side of  the estimate (\ref{estim-base-NS}) we take the limit when $R\to +\infty$ to obtain the identity $\U=0$. Moreover,  by the identities  $\ds{q=\sum_{i,j=1}^{3} \mathcal{R}_i \mathcal{R}_j (u_i\, u_j)}$ and $\ds{\nabla q=\nabla p}$, we conclude that $p$ is a constant.
\item[$2)$] The case  	$\eta(\gamma,p)=0$. In this case, by the estimate (\ref{estim-base-NS})  we have 
\[  \int_{B_{R/2}} \vert \nabla \otimes \U \vert^2 dx \leq  \,\, c\, \Vert \U \Vert^{2}_{M^{p}_{\gamma}}+ c\,\Vert \U \Vert^{3}_{M^{p}_{\gamma}},
\]
\end{enumerate}	
and taking the limit when $R\to +\infty$ we obtain $\ds{\int_{\Rt} \vert \nabla \otimes \U \vert^2 dx \leq  \,\, c\, \Vert \U \Vert^{2}_{M^{p}_{\gamma}}+ c\,\Vert \U \Vert^{3}_{M^{p}_{\gamma}}}$, hence we can write $\U \in \dot{H}^{1}(\Rt)$. We will use now the additional hypothesis $\U \in \dot{B}^{-1}_{\infty,\infty}(\Rt)$ to conclude the identity $\U=0$. 

\medskip

Indeed, with the information  $\U \in \dot{B}^{-1}_{\infty, \infty}(\Rt)$ we can apply the improved Sobolev inequalities (see the article \cite{GerardMeyerOru} for a proof of these inequalities) and we write $\ds{\Vert \U \Vert_{L^4} \leq c \Vert \U \Vert^{\frac{1}{2}}_{\dot{H}^{1}} \Vert \U \Vert^{\frac{1}{2} }_{\dot{B}^{-1}_{\infty,\infty}}}$. Once we dispose of the information $\U\in L^4(\Rt)$ we can derive now the identity $\U=0$ as follows: multiplying equation (\ref{NS})  by $\U$ and integrating on the whole space $\Rt$ we have 
$$ \int_{\Rt}(-\Delta \U)\cdot \U dx= \int_{\Rt} ((\U \cdot \nabla) \U )\cdot \U dx + \int_{\Rt} \nabla p  \cdot \U dx,$$ where due to the fact  $\U \in \dot{H}^{1}\cap L^4(\Rt)$ each term in this identity is well-defined.  Indeed, for the term in the left-hand side remark that as $\U \in \dot{H}^{1}(\Rt)$ then  we have $-\Delta \U \in \dot{H}^{-1}(\Rt)$.  Then, for the first term in the right-hand side, as $div(\U)=0$ we write $(\U \cdot \vec{\nabla}) \U=div(\U \otimes \U)$ where, as $\U \in L^{4}(\Rt)$ by the H\"older inequalities we have $\U \otimes \U \in L^{2}(\Rt)$ and then $div(\U \otimes \U)\in \dot{H}^{-1}(\Rt)$. Finally, in order to study the second term in the right-hand side, we write the pressure $p$ as $p= \frac{1}{-\Delta} div (div (\U \otimes \U))$ hence we get $p \in L^{2}(\Rt)$ (since we have $\U \otimes \U \in L^{2}(\Rt)$) and then $\nabla p \in \dot{H}^{-1}(\Rt)$.

\medskip 

Now, integrating by parts each term in the identity above we have that $\int_{\Rt}(-\Delta \U)\cdot \U dx = \int_{\Rt}\vert \nabla \otimes \U \vert^2 dx$, and moreover $\int_{\Rt} ((\U \cdot \nabla) \U )\cdot \U dx=0$ and $\int_{\Rt} \nabla p  \cdot \U dx=0$. With these identities we get $\int_{\Rt}\vert \nabla \otimes \U \vert^2 dx =0$ and thus we have $\U=0$.  Proposition \ref{Prop:NS} is now proven. \finpv

\section{The non-stationary case}\label{Sec:Non-Stat}

\subsection{Proof of Theorem \ref{Th2}}

Recall first that by hypothesis of Theorem \ref{Th2} we have $\vu \in M^{p}_{\gamma,0}L^p(0,T)$ and $\nabla \otimes \v \in  M^{p}_{\gamma}L^p(0,T)$, where the parameters $0<\gamma<3\leq p <+\infty$ are such that $\eta(\gamma,p)\leq 0$.  Then, by Proposition \ref{Prop-presion} we have the identity $\nabla p = \nabla q$, where the quantity $q$ in defined in expression (\ref{Presion}); and from now on we will consider the equations (\ref{EickLes-non-stat}) with the term $\nabla q$ instead of the term $\nabla p$. 

\medskip 
 
We will apply the local energy balance  (\ref{energ-loc}) to a suitable test function and for this we will follow some of the ideas of \cite{PFPG1}. Let $0<t_0<t_1<T$. For a parameter $\varepsilon>0$, we will consider a function $\alpha_{\varepsilon,t_0,t_1}(t)$ which converges \emph{a.e.} to $\ds{\mathds{1}_{[t_0, t_1]}(t)}$ and such that $\ds{\frac{d}{dt} \alpha_{\varepsilon,t_0,t_1}(t)}$ is the difference between two identity approximations: the first one in $t_0$ and the second one in $t_1$.  For this, let $\alpha \in \mathcal{C}^{\infty}(\mathbb{R})$ be a function such that $\alpha(t)=0$ for $-\infty < t <1/2$ and $\alpha(t)=1$ for $1 < t <+\infty$. Then, for $\varepsilon< \min (t_0 / 2, T-t_1)$ we set the function $\ds{\alpha_{\varepsilon,t_0,t_1} (t)= \alpha \left( \frac{t-t_0}{\varepsilon}\right)- \alpha \left( \frac{t-t_1}{\varepsilon}\right)}$.\\

On the other hand, for  $R\geq 1$ let $\theta_R(x)$ be function test given in (\ref{funTheta}).Then, we consider the function test $\ds{\alpha_{\varepsilon,t_0,t_1}(t)\theta_R(x)}$ and by (\ref{energ-loc}) we can write 
\begin{equation*}
\begin{split}
&- \int_{\R}\int_{\Rt} \frac{ \vert \vu \vert^2 + \vert \nabla \otimes \v \vert^2}{2} \partial_{s} \alpha_{\varepsilon,t_0,t_1} \theta_R \, dx\, ds + \int_{\R} \int_{\Rt} \vert \nabla \otimes \vu \vert^2 \alpha_{\varepsilon,t_0,t_1} \theta_R dx\, ds  + \int_{\R}\int_{\Rt} \vert \Delta \v \vert^2 \alpha_{\varepsilon,t_0,t_1} \theta_R dx\, ds \\
\leq & \int_{\R}\int_{\Rt}  \left( \frac{\vert \vu \vert^2+ \vert \nabla \otimes \v \vert^2}{2} \right) \alpha_{\varepsilon,t_0,t_1}  \Delta\theta_R dx \, ds +  \int_{\R}\int_{\Rt}  \left( \left[ \frac{\vert \vu \vert^2 +\vert \nabla \otimes \v \vert^2}{2}+q \right] \vu  \right) \cdot \alpha_{\varepsilon,t_0,t_1} \nabla \theta_R dx \, ds \\
&  \int_{\R} \int_{\Rt} \sum_{k=1}^{3}  ([\vu \cdot \nabla)\, \v ] \cdot \partial_k \v ) \alpha_{\varepsilon,t_0,t_1} \partial_k\theta_R dx\, ds  -  \int_{\R} \int_{\Rt}  \vert \nabla \otimes \v \vert^2 \v \cdot \Delta \v \, \alpha_{\varepsilon,t_0,t_1}\theta_R dx \, ds.
\end{split}
\end{equation*}
Now, taking the limit when $\varepsilon \to 0$, by the dominated convergence theorem we obtain (when the limit in the left side is well-defined)
\begin{equation*}
\begin{split}
&- \lim_{\varepsilon\to 0 }\int_{\R}\int_{\Rt} \frac{ \vert \vu \vert^2 + \vert \nabla \otimes \v \vert^2}{2} \partial_{s} \alpha_{\varepsilon,t_0,t_1} \theta_R \, dx\, ds + \int_{t_0}^{t_1} \int_{\Rt} \vert \nabla \otimes \vu \vert^2  \theta_R dx\, ds  + \int_{t_0}^{t_1}\int_{\Rt} \vert \Delta \v \vert^2  \theta_R dx\, ds \\
\leq & \int_{t_0}^{t_1}\int_{\Rt}  \left( \frac{\vert \vu \vert^2+ \vert \nabla \otimes \v \vert^2}{2} \right)  \Delta\theta_R dx \, ds +  \int_{t_0}^{t_1}\int_{\Rt}  \left( \left[ \frac{\vert \vu \vert^2 +\vert \nabla \otimes \v \vert^2}{2}+q \right] \vu  \right) \cdot  \nabla \theta_R dx \, ds \\
&  \int_{t_0}^{t_1} \int_{\Rt} \sum_{k=1}^{3}  ([\vu \cdot \nabla)\, \v ] \cdot \partial_k \v )  \partial_k\theta_R dx\, ds  -  \int_{t_0}^{t_1} \int_{\Rt}  \vert \nabla \otimes \v \vert^2 \v \cdot \Delta \v \,\theta_R dx \, ds.
\end{split}
\end{equation*}
At this point, we must study the expression $\ds{- \lim_{\varepsilon\to 0 }\int_{\R}\int_{\Rt} \frac{ \vert \vu \vert^2 + \vert \nabla \otimes \v \vert^2}{2} \partial_{s} \alpha_{\varepsilon,t_0,t_1} \theta_R \, dx\, ds}$. To make the writing more simple, let us define the function 
$\ds{A_R(s)= \int_{\Rt} \frac{ \vert \vu(s,x) \vert^2 + \vert \nabla \otimes \v(s,x) \vert^2}{2}  \theta_R \, dx}$. Then, assuming that $t_0$ and $t_1$ are Lebesgue points of the function $A_R(s)$, and moreover, since  
$$ \int_{\R}\int_{\Rt} \frac{ \vert \vu \vert^2 + \vert \nabla \otimes \v \vert^2}{2} \partial_{s} \alpha_{\varepsilon,t_0,t_1} \theta_R \, dx\, ds = -\frac{1}{2} \int_{\R} A_R(s) \partial_{s} \alpha_{\varepsilon,t_0,t_1} ds,$$
then we have 
$$ - \lim_{\varepsilon\to 0 }\int_{\R}\int_{\Rt} \frac{ \vert \vu \vert^2 + \vert \nabla \otimes \v \vert^2}{2} \partial_{s} \alpha_{\varepsilon,t_0,t_1} \theta_R \, dx\, ds= \frac{1}{2}(A_R(t_1)- A_R(t_0)).$$
On the other hand, recall that by point $4)$ in Definition \ref{def-suitable} we have that the functions $\vu(t,\cdot)$ and $\nabla \otimes \v(t,\cdot)$ are strong continuous at $t=0$  and then we can replace $t_0$ by $0$. Moreover, for $0<t<T$, always by point $4)$ in Definition \ref{def-suitable}  we have that the functions $\vu(t,\cdot)$ and $\nabla \otimes \v(t,\cdot)$ are weak continuous at $t$ and then we obtain 
$\ds{A_R(t)\leq \liminf_{t_1\to t} A_R(t_1)}$. Thus, we can also replace $t_1$ for $t$. \\

With this information, for every $0\leq t\leq T$ we can write 
\begin{equation*}
\begin{split}
&\int_{\Rt} \frac{ \vert \vu(t,\cdot) \vert^2 + \vert \nabla \otimes \v(t,\cdot) \vert^2}{2}  \theta_R \, dx + \int_{0}^{t} \int_{\Rt} \vert \nabla \otimes \vu \vert^2  \theta_R dx\, ds  + \int_{0}^{t}\int_{\Rt} \vert \Delta \v \vert^2  \theta_R dx\, ds \\
\leq & \int_{\Rt} \frac{ \vert \vu_0 \vert^2 + \vert \nabla \otimes \v_0 \vert^2}{2}  \theta_R \, dx +  \int_{0}^{t}\int_{\Rt}  \left( \frac{\vert \vu \vert^2+ \vert \nabla \otimes \v \vert^2}{2} \right)  \Delta\theta_R dx \, ds\\
& +  \int_{0}^{t}\int_{\Rt}  \left( \left[ \frac{\vert \vu \vert^2 +\vert \nabla \otimes \v \vert^2}{2}+q \right] \vu  \right) \cdot  \nabla \theta_R dx \, ds + \int_{0}^{t} \int_{\Rt} \sum_{k=1}^{3}  ([\vu \cdot \nabla)\, \v ] \cdot \partial_k \v )  \partial_k\theta_R dx\, ds \\
& -  \int_{0}^{t} \int_{\Rt}  \vert \nabla \otimes \v \vert^2 \v \cdot \Delta \v \,\theta_R dx \, ds.
\end{split}
\end{equation*}

In this inequality we must  study now  the term $\ds{-  \int_{0}^{t} \int_{\Rt}  \vert \nabla \otimes \v \vert^2 \v \cdot \Delta \v \theta_R dx \, ds}$. Recall  that by (\ref{Iden}) we have the identity (in the distributional sense) $\ds{\vert \nabla \otimes \v \vert^2=-\v \cdot \Delta \v}$, moreover,  as we have $\ds{\vert \v \vert=1}$,  then we can write 
$$ -  \int_{0}^{t} \int_{\Rt}  \vert \nabla \otimes \v \vert^2 \v \cdot \Delta \v \theta_R dx \, ds =  \int_{0}^{t} \int_{\Rt}  \vert \v \cdot \Delta \v \vert^2 \theta_R dx \, ds \leq \int_{0}^{t} \int_{\Rt}  \vert \v \vert^2 \vert  \Delta \v \vert^2 \theta_R dx \, ds  \leq \int_{0}^{t} \int_{\Rt}  \vert  \Delta \v \vert^2 \theta_R dx \, ds.$$
By this estimate and the previous inequality we get 
 \begin{equation*} 
 \begin{split}
 &\int_{\Rt} \frac{ \vert \vu \vert^2 + \vert \nabla \otimes \v \vert^2}{2}\theta_R \, dx + \int_{0}^{t} \int_{\Rt} \vert \nabla \otimes \vu \vert^2 \theta_R dx\, ds  \leq  \int_{\Rt} \frac{ \vert \vu_0 \vert^2 + \vert \nabla \otimes \v_0 \vert^2}{2}\theta_R \, dx  \\
 &+ \int_{0}^{t}\int_{\Rt}\left( \frac{\vert \vu \vert^2+ \vert \nabla \otimes \v \vert^2}{2} \right)  \Delta \theta_R dx \, ds +\int_{0}^{t}\int_{\Rt}  \left( \left[ \frac{\vert \vu \vert^2 +\vert \nabla \otimes \v \vert^2}{2}+q \right] \vu  \right) \cdot \nabla \theta_R dx \, ds \\
 & - \int_{0}^{s} \int_{\Rt} \sum_{k=1}^{3} \partial_k ([\vu \cdot \nabla)\, \v ] \cdot \partial_k \v ) \theta_R dx\, ds.
 \end{split}
 \end{equation*}
Now, as we  have $u_0 \in L^2(\Rt)$ and $\v_0 \in \dot{H}^{1}(\Rt)$, and moreover,  recalling that  $supp(\theta_R)\subset B_R$, $supp( \nabla \theta_R)\subset \mathcal{C}(R/2, R)$ and $supp( \Delta  \theta_R)\subset \mathcal{C}(R/2, R)$, then  we write  
\begin{eqnarray}\label{Desigualdad} \nonumber
& & \int_{B_R} \frac{ \vert \vu \vert^2 + \vert \nabla \otimes \v \vert^2}{2}\theta_R \, dx + \int_{0}^{t} \int_{B_R} \vert \nabla \otimes \vu \vert^2 \theta_R dx\, ds   \leq  \Vert u_0 \Vert^{2}_{L^2}+\Vert \v_0 \Vert^{2}_{\dot{H}^1} \\ \nonumber
& &  + \int_{0}^{t}\int_{\mathcal{C}(R/2,R)}  \left( \frac{\vert \vu \vert^2+ \vert \nabla \otimes \v \vert^2}{2} \right) \Delta\theta_R dx \, ds  + \int_{0}^{t}\int_{\mathcal{C}(R/2,R)} \left( \left[ \frac{\vert \vu \vert^2 +\vert \nabla \otimes \v \vert^2}{2}+q \right] \vu  \right) \nabla \theta_R dx \, ds \\  \nonumber
& & + \int_{0}^{t} \int_{\mathcal{C}(R/2,R)} \sum_{k=1}^{3}  ([\vu \cdot \nabla)\, \v ] \cdot \partial_k \v ) \partial_k \theta_R dx\, ds \\
&=& \Vert u_0 \Vert^{2}_{L^2}+\Vert \v_0 \Vert^{2}_{\dot{H}^1}+ I_1+I_2+I_3, 
\end{eqnarray} where we will show that we have $\ds{\lim_{R \to +\infty} I_i=0}$ for $i=1,2,3$. Indeed, for the term $I_1$  recall that we have $\ds{\Vert \Delta \theta_R \Vert_{L^{\infty}} \leq \frac{c}{R^2}}$, and the we get 
$$ I_1 \leq \frac{c}{R^2} \int_{0}^{t} \int_{\mathcal{C}(R/2, R)} (\vert \vu \vert^2 + \vert \nabla \otimes \v \vert^2) dx\, ds \leq c\, R^{1-6/p} \int_{0}^{t} \left( \int_{\mathcal{C}(R/2, R)} (\vert \vu \vert^p + \vert \nabla \otimes \v \vert^p)  dx \right)^{2/p} \, ds,$$ thereafter, by the H\"older inequalities in the temporal variable (with $1=2/p+(p-2)/p$), and moreover, recalling that we define $\eta(\gamma,p)=\gamma/p-3/p+2/3$,  we have  

\begin{equation*}
\begin{split}
& c\, R^{1-6/p} \int_{0}^{t} \left( \int_{\mathcal{C}(R/2, R)} (\vert \vu \vert^p + \vert \nabla \otimes \v \vert^p)  dx \right)^{2/p} \, ds\\
\leq &  c\, R^{1-6/p} \left(  \int_{0}^{t} \int_{\mathcal{C}(R/2, R)} (\vert \vu \vert^p + \vert \nabla \otimes \v \vert^p)  dx\, ds \right)^{2/p} t^{(p-2)/p} \\
\leq &  c\, R^{1-6/p+2\gamma/p} \left(  \frac{1}{R^{\gamma}} \int_{0}^{t} \int_{\mathcal{C}(R/2, R)} (\vert \vu \vert^p + \vert \nabla \otimes \v \vert^p)  dx\, ds \right)^{2/p} t^{(p-2)/p}\\
\leq & c\, R^{2(1/2-3/p+\gamma/p)}  \left(  \frac{1}{R^{\gamma}} \int_{0}^{t} \int_{\mathcal{C}(R/2, R)} (\vert \vu \vert^p + \vert \nabla \otimes \v \vert^p)  dx\, ds \right)^{2/p} t^{(p-2)/p} \\
\leq & c\, R^{2(2/3-3/p+\gamma/p -1/6)}  \left(  \frac{1}{R^{\gamma}} \int_{0}^{t} \int_{\mathcal{C}(R/2, R)} (\vert \vu \vert^p + \vert \nabla \otimes \v \vert^p)  dx\, ds \right)^{2/p} t^{(p-2)/p}\\
\leq & c\, R^{2\eta(\gamma,p)-1/3} \left(  \frac{1}{R^{\gamma}} \int_{0}^{t} \int_{\mathcal{C}(R/2, R)} (\vert \vu \vert^p + \vert \nabla \otimes \v \vert^p)  dx\, ds \right)^{2/p} T^{(p-2)/p}. 
\end{split}
\end{equation*}
Hence, as $\eta(\gamma,p)\leq 0$ and $R\geq 1$,  we finally write 
$$ I_1 \leq  c \, \frac{T^{(p-2)/p}}{R^{1/3}}  \Vert \U \Vert^{2}_{M^{p}_{\gamma}L^{p}(0,T)} + c \frac{T^{(p-2)/p}}{R^{1/3}}   \Vert \nabla \otimes \v \Vert^{2}_{M^{3}_{1}L^{3}(0,T)}.$$
But, as we have the information $\vu \in M^{p}_{\gamma,0}L^{p}(0,T)$ and $\nabla \otimes \v \in M^{p}_{\gamma}L^p(0,T)$,  taking the limit when $R\to +\infty$   we obtain 
$\ds{\lim_{R\to +\infty}I_1 =0}$. \\ 

For the term $I_2$, by the  estimates (\ref{I1}), (\ref{I2}) and (\ref{I3}), we have 
\begin{eqnarray*}
I_2 &\leq & \int_{0}^{t} \left( \int_{\mathcal{C}(R/2,R)} \vert \vu \vert^p dx \right)^{2/p} R^{2-9/p} \left( \int_{\mathcal{C}(R/2,R)} \vert \vu \vert^p dx \right)^{1/p} ds \\
& & + \int_{0}^{t} \left( \int_{\mathcal{C}(R/2,R)} \vert \nabla \otimes \v \vert^p dx \right)^{2/p} R^{2-9/p}  \left( \int_{\mathcal{C}(R/2,R)} \vert \vu \vert^p dx \right)^{1/p} ds \\
& & + \int_{0}^{t} \left( \int_{\mathcal{C}(R/2,R)} \vert q \vert^{p/2} dx \right)^{2/p} R^{2-9/p} \left( \int_{\mathcal{C}(R/2,R)} \vert \vu \vert^3 dx \right)^{1/p} ds,
\end{eqnarray*} 
hence, since $\eta(\gamma,p)=\gamma/p-3/p+2/3$ then  we write 
\begin{eqnarray*}
	I_2 &\leq & \int_{0}^{t} \left(\frac{1}{R^\gamma} \int_{\mathcal{C}(R/2,R)} \vert \vu \vert^p dx \right)^{2/p} R^{2-9/p+3\gamma/p} \left( \frac{1}{R^\gamma} \int_{\mathcal{C}(R/2,R)} \vert \vu \vert^p dx \right)^{1/p} ds \\
	& & + \int_{0}^{t} \left( \frac{1}{R^\gamma}  \int_{\mathcal{C}(R/2,R)} \vert \nabla \otimes \v \vert^p dx \right)^{2/p} R^{2-9/p+3\gamma/p}  \left(\frac{1}{R^\gamma}  \int_{\mathcal{C}(R/2,R)} \vert \vu \vert^p dx \right)^{1/p} ds \\
	& & + \int_{0}^{t} \left( \frac{1}{R^\gamma}  \int_{\mathcal{C}(R/2,R)} \vert q \vert^{p/2} dx \right)^{2/p} R^{2-9/p+3\gamma/p} \left( \frac{1}{R^\gamma}  \int_{\mathcal{C}(R/2,R)} \vert \vu \vert^3 dx \right)^{1/p} ds \\
	&\leq & R^{3\eta(\gamma,p)} \left[  \int_{0}^{t} \left(\frac{1}{R^\gamma} \int_{\mathcal{C}(R/2,R)} \vert \vu \vert^p dx \right)^{2/p}   \left( \frac{1}{R^\gamma} \int_{\mathcal{C}(R/2,R)} \vert \vu \vert^p dx \right)^{1/p} ds \right. \\
	& & +   \int_{0}^{t} \left( \frac{1}{R^\gamma}  \int_{\mathcal{C}(R/2,R)} \vert \nabla \otimes \v \vert^p dx \right)^{2/p}  \left(\frac{1}{R^\gamma}  \int_{\mathcal{C}(R/2,R)} \vert \vu \vert^p dx \right)^{1/p} ds \\
	& &  \left. \int_{0}^{t} \left( \frac{1}{R^\gamma}  \int_{\mathcal{C}(R/2,R)} \vert q \vert^{p/2} dx \right)^{2/p} \left( \frac{1}{R^\gamma}  \int_{\mathcal{C}(R/2,R)} \vert \vu \vert^3 dx \right)^{1/p} ds \right].
\end{eqnarray*} 

Using first  the fact that $\eta(\gamma,p)\leq 0$, and moreover,  applying the H\"older inequalities in the temporal variable (with $1=2/p+1/p+(p-3)/p$), we obtain
\begin{eqnarray*}
	I_2 &\leq & \left( \frac{1}{R^\gamma} \int_{0}^{t}  \int_{\mathcal{C}(R/2,R)} \vert \vu \vert^p dx\, ds \right)^{2/p}  \left(\frac{1}{R^\gamma} \int_{0}^{t} \int_{\mathcal{C}(R/2,R)} \vert \vu \vert^p dx \, ds\right)^{1/p} t^{(p-3)/p}  \\
	& & + \left( \frac{1}{R^\gamma} \int_{0}^{t}  \int_{\mathcal{C}(R/2,R)} \vert \nabla \otimes \v \vert^p dx \, ds\right)^{2/p}  \left(\frac{1}{R^\gamma} \int_{0}^{t} \int_{\mathcal{C}(R/2,R)} \vert \vu \vert^p dx \, ds\right)^{1/p} t^{(p-3)/p} \\
	& & + \left( \frac{1}{R^\gamma} \int_{0}^{t}  \int_{\mathcal{C}(R/2,R)} \vert q \vert^{p/2} dx \, ds\right)^{2/p}  \left( \frac{1}{R^\gamma}\int_{0}^{t}\int_{\mathcal{C}(R/2,R)} \vert \vu \vert^p dx\, ds \right)^{1/p} t^{(p-3)/p}\\
&\leq & \left( \Vert \vu \Vert^{2}_{M^{p}_{\gamma}L^p(0,T)}   + \Vert \nabla \otimes \v  \Vert^{2}_{M^{p}_{\gamma}L^p(0,T)} +  \Vert q \Vert_{M^{p/2}_{\gamma}L^{p/2}(0,T)}\right) \left(\frac{1}{R^\gamma} \int_{0}^{T } \int_{\mathcal{C}(R/2,R)} \vert \vu \vert^p dx \, ds\right)^{1/p} \, T^{(p-3)/p}. 
\end{eqnarray*} 
At this point, as we have  $\vu \in M^{p}_{\gamma,0}L^{p}(0,T)$ and $\nabla \otimes \v \in M^{p}_{\gamma}L^p(0,T)$, the  by point $1)$ of Lemma \ref{Lema-Tech2} we get    
\begin{equation*}\label{p-non-stat}
\Vert q \Vert_{M^{p/2}_{\gamma}L^{p/2}(0,T)} \leq c \left(\Vert \vu \Vert^{2}_{M^{p}_{\gamma}L^{p}(0,T)}+\Vert \nabla \otimes \v \Vert^{2}_{M^{p}_{\gamma}L^{p}(0,T)}\right).
\end{equation*} 
Thus, getting back to the previous estimate  we can write 
\begin{eqnarray*}
	I_2 &\leq & c \left( \Vert \vu \Vert^{2}_{M^{p}_{\gamma}L^p(0,T)}   + \Vert \nabla \otimes \v  \Vert^{2}_{M^{p}_{\gamma}L^p(0,T)} \right) \left( \frac{1}{R^\gamma}\int_{0}^{T } \int_{\mathcal{C}(R/2,R)} \vert \vu \vert^p dx \, ds\right)^{1/p}\, T^{(p-3)/p},  
\end{eqnarray*} and then, taking the limit when $R\to +\infty$
we have $\ds{\lim_{R\to +\infty}I_2 =0}$. 

\medskip 

Finally, for the term $I_3$, applying the H\"older inequalities in the spatial variable (with $1= 1/p+2/p+(p-3)/p$), we have 
\begin{equation*}
\begin{split}
I_3  & = \sum_{i,j,k=1}^{3} \int_{0}^{t}\int_{\mathcal{C}(R/2,R)} (u_j \partial_j v_i) (\partial_k v_i) \partial_k \theta_R  dx\,ds \leq  c \int_{0}^{t} \int_{\mathbb{C}(R/2, R)} \vert \vu \vert \vert \nabla \otimes \v \vert^2 \vert \nabla \theta_R \vert dx\,ds \\
&\leq c\, \int_{0}^{t} \left( \int_{\mathcal{C}(R/2,R) } \vert \vu \vert^p dx \right)^{1/p}  \left( \int_{\mathcal{C}(R/2,R) } \vert \nabla \otimes  \v \vert^p dx \right)^{2/p} \left( \int_{\mathcal{C}(R/2,R) } \vert \nabla \theta_R \vert^{p/(p-3)} dx \right)^{(p-3)/p}\, ds.   
\end{split}
\end{equation*} 
Moreover, in the last term, as $\Vert \nabla \theta_R \Vert_{L^{\infty}} \leq \frac{c}{R}$  then we can write 
\begin{equation*}
\begin{split}
I_3 & \leq c\, \left[  \int_{0}^{t} \left( \int_{\mathcal{C}(R/2,R) } \vert \vu \vert^p dx \right)^{1/p}  \left( \int_{\mathcal{C}(R/2,R) } \vert \nabla \otimes  \v \vert^p dx \right)^{2/p} ds \right]\,  R^{2-9/p}\\
&\leq  c\, \left[  \int_{0}^{t} \left( \frac{1}{R^\gamma} \int_{\mathcal{C}(R/2,R) } \vert \vu \vert^p dx \right)^{1/p}  \left( \frac{1}{R^\gamma}  \int_{\mathcal{C}(R/2,R) } \vert \nabla \otimes  \v \vert^p dx \right)^{2/p} ds \right]\,  R^{2-9/p+3\gamma/p}.
\end{split}
\end{equation*}
Then, recalling that $\eta(\gamma,p)=2/3-3/p+\gamma/p$, and moreover, as we assume $\eta(\gamma,p)\leq 0$,  we obtain 
\begin{equation*}
\begin{split}
I_3 & \leq   c\,  \left[\int_{0}^{t} \left( \frac{1}{R^\gamma} \int_{\mathcal{C}(R/2,R) } \vert \vu \vert^p dx \right)^{1/p}  \left( \frac{1}{R^\gamma}  \int_{\mathcal{C}(R/2,R) } \vert \nabla \otimes  \v \vert^p dx \right)^{2/p} ds \right]\,  R^{3\eta(\gamma,p)}\\
	& \leq c\,  \int_{0}^{t} \left( \frac{1}{R^\gamma} \int_{\mathcal{C}(R/2,R) } \vert \vu \vert^p dx \right)^{1/p}  \left( \frac{1}{R^\gamma}  \int_{\mathcal{C}(R/2,R) } \vert \nabla \otimes  \v \vert^p dx \right)^{2/p}ds.
\end{split}	
\end{equation*}
We apply now the H\"older inequalities in the temporal variable (with $1=1/p+2/p+(p-3)/p$) to write 
\begin{equation*}
\begin{split}
I_3 & \leq   c\, \left(   \frac{1}{R^\gamma}  \int_{0}^{t}  \int_{\mathcal{C}(R/2,R) } \vert \vu \vert^p dx\, ds \right)^{1/p}  \left(  \frac{1}{R^\gamma}  \int_{0}^{t}   \int_{\mathcal{C}(R/2,R) } \vert \nabla \otimes  \v \vert^p dx\, ds  \right)^{2/p}\, t^{(p-3)/p} \\
& \leq c\,  \left(   \frac{1}{R^\gamma}  \int_{0}^{t}  \int_{\mathcal{C}(R/2,R) } \vert \vu \vert^p dx\, ds \right)^{1/p}  \, \Vert \nabla \otimes \v \Vert^{2}_{M^{p}_{\gamma}L^p(0,T)} \, T^{(p-3)/p}. 
\end{split}	
\end{equation*}
Hence, as $\vu \in M^{p}_{\gamma,0}L^p(0,T)$,  taking the limit when $R\to +\infty$ we obtain $\ds{\lim_{R\to +\infty} I_3=0}$. \\

Once we have proven that $\ds{\lim_{R\to +\infty} I_i=0}$ for $i=1,2,3$, we get back to (\ref{Desigualdad}) where we take the limit when $R\to +\infty$,  and thus for $0\leq t \leq T$ we get the global energy inequality (\ref{global-energ}). Theorem \ref{Th2} is proven. \finpv
\subsection{Proof of Corollary \ref{Col1}} This proof is straightforward. We just observe that by the global energy inequality (\ref{global-energ}) if the initial datum verify $\vu_0=0$ and $\nabla \otimes \v_0=0$ then for all time $0<t\leq T$ we have 
$\Vert \vu(t,\cdot )\Vert^{2}_{L^2}=0$ and $\Vert \v(t,\cdot)\Vert^{2}_{\dot{H}^{1}}=0$, hence $\vu=0$ and $\nabla \otimes \v=0$ on $[0,T]\times \Rt$.  \finpv

	\begin{appendices}
	\section{Appendix}\label{AppendixA}
	Consider the velocity field $\U$, the pressure term $p$ and the vector field $\V$ defined as follows: 
	\begin{equation}\label{ExStat}
	\begin{split}
	\U(x_1, x_2, & x_3)= (2x_1, 2x_2, -4x_3),\quad p(x_1, x_2, x_3)= -(2x^{2}_{1} +2x^{2}_{2}+8 x^{3}_{2}), \\ 
	& \text{and}\quad \V(x_1, x_2, x_3)=\left\{\begin{array}{ll}\vspace{2mm} 
	(x_1,  x_2, 0), \,\,\, \text{if}\,\,\, x^{2}_{1}+x^{2}_{2}=1, \\
	\left(\frac{1}{\sqrt{3}}, \frac{1}{\sqrt{3}}, \frac{1}{\sqrt{3}}\right),  \,\,\, \text{if}\,\,\, x^{2}_{1}+x^{2}_{2}\neq 1. 
	\end{array}\right.
	\end{split}
	\end{equation}
	We have $\vert \V \vert=1$, and  using some basic rules of the vector calculus we easily get  that the triplet $(\U, p, \V)$ defined as above is also  a solution  of the system (\ref{EickLes}).   Indeed, for the first equation in the system (\ref{EickLes}), by definition of the  vector field $\v(x_1, x_2, x_2)$ above we observe first that we have 
	$\ds{ \text{div} ( \nabla \otimes  \V \odot \nabla \otimes \V)=0}$. 
	Thereafter, if we set the scalar field $\psi(x_1,x_2,x_3)=x^{2}_{1}+x^{2}_{2}-2x^{2}_{3}$ we may observe that we have $\U = \nabla \psi$ and moreover we have $\Delta \psi =0$. With these identities we can write the following computations. First, we have $\ds{\Delta \U = \Delta (\nabla \psi)= \nabla (\Delta \psi)=0}$. On the other hand we have 
	$\ds{ (\U \cdot \nabla) \U = \frac{1}{2} \vert \U \vert + (\nabla \wedge \U) \wedge \U}$, and as $\U = \nabla \psi$ then we get 
	$\ds{ (\U \cdot \nabla) \U = \frac{1}{2} \vert \U \vert}$. 
	With this identity and the definition of the pressure term $p$ given above we find that $(\U,p,\V)$ verify the first equation of (\ref{EickLes}). Moreover we have $div(\U)=0$. \\
	
	For the second equation in (\ref{EickLes}), we observe first that for the case $x^{2}_{1}+x^{2}_{2}\neq 1$ the vector field $\v(x_1,x_2,x_3)$ defined above is a constant vector and then the second equation in (\ref{EickLes}) trivially holds. For the other case, when $x^{2}_{1}+x^{2}_{2} = 1$,  we have  $\v(x_1, x_2, x_3)=(x_1,x_2,0)$ and then we get $\Delta \v=0$. Moreover, by definition of the vector field $\U(x_1,x_2,x_3)$ is easy to see that we have the identity $(\U \cdot \nabla) \V= \vert \nabla \otimes \V\vert^2 \V$. Then, the second equation in (\ref{EickLes}) also  holds true. 
	
\section{Appendix}\label{AppendixB}
Here we give a proof of the embedding (\ref{embedding}).  It is enough to prove the last inclusion in this embedding, and for this,  for  all  $R\geq 1$  and $t \geq 0$ fixed,  we have the estimate 
\[ \int_{B_R} \vert f(t,x) \vert^p  dx \leq c R^{3(1-p/r)} \Vert f(t,\cdot)\Vert^{p}_{L^{r,\infty}} \leq c R^{3(1-p/r)} \Vert f(t,\cdot)\Vert^{p}_{L^{r,q}},\]
with $r<q<+\infty$. For a proof of this estimate the Proposition $1.1.10$, page 21 of the book \cite{DCh}. Then,  for $0<\gamma<3$ we write 
\[  \frac{1}{R^{\gamma}}\int_{B_R} \vert f(t,x) \vert^p  dx \leq c R^{3(1-p/r)-\gamma} \Vert f(t,\cdot)\Vert^{p}_{L^{r,q}}. \]
In the right side of this estimate we impose now the condition $\ds{3(1-p/r)-\gamma \leq 0}$, which is equivalent to the inequality $\ds{2/3-3/r \leq \gamma/p -3/p +2/3=\eta(\gamma,p)}$. Moreover, as we assume $\eta(\gamma,p) \leq 0$ we get $2/3-3/r \leq 0$ which give us the restriction on the parameter $r$: $r\leq 9/2$. 

\medskip

Thus, when $\ds{3(1-p/r)-\gamma \leq 0}$ holds, we can write 
\[  \frac{1}{R^{\gamma}}\int_{B_R} \vert f(t,x) \vert^p  dx \leq c R^{3(1-p/r)-\gamma} \Vert f(t,\cdot)\Vert^{p}_{L^{r,\infty}}\leq c  \Vert f(t,\cdot)\Vert^{p}_{L^{r,q}},\] hence, integrating in the interval of time $[0,T]$, and moreover,  following the same to prove (\ref{Emb}),  we finally obtain the   last inclusion in  (\ref{embedding}). 
\end{appendices}

\end{document}